\documentclass[12pt]{svmult}
\usepackage{amssymb}
\usepackage{amsmath}
\usepackage{lineno}

\usepackage{graphicx,multicol,xcolor}
\usepackage{amsmath, amscd, amssymb, amsfonts, dsfont}
\usepackage{pstricks,pst-plot}

\usepackage{mathptmx}

\makeatletter

\newenvironment{figurehere}
  {\def\@captype{figure}}
  {}
\makeatother

\usepackage{pstricks,pst-plot}
\usepackage{graphicx,multicol,xcolor}

\usepackage{mathptmx}
\usepackage{helvet}
\usepackage{courier}
\usepackage{type1cm}

\usepackage{makeidx}

\usepackage[bottom]{footmisc}

\usepackage{comment}

\makeindex             % used for the subject index
\newcommand{\crit}[1]{\textup{Crit}(#1)}
\newcommand{\grad}[1]{\textup{grad}(#1)}
\newcommand{\ind}[2]{\textup{ind}(#1;#2)}

\begin{document}

%%%%%%%%%%%%%%%%%%%%%%%%%%
%%%%%%%%%%%%%%%%%%%%%%%%%%
%\title{A Quick View of Lagrangian Floer Homology }
%\author{Andr\'es Pedroza}
%\address{Facultad de Ciencias\\
 %          Universidad de Colima\\  
%           Bernal D\'{\i}az del Castillo No. 340\\
%           Colima, Col., Mexico 28045}
%\email{andres\_pedroza@ucol.mx}

%%%%%%%%%%%%%%%%%%%%%%%%%%
%%%%%%%%%%%%%%%%%%%%%%%%%%

\title*{A Quick View of Lagrangian Floer Homology}
\author{Andr\'es Pedroza}
\institute{Andr\'es Pedroza \at 
Facultad de Ciencias, 
Universidad de Colima, 
Bernal D\'{\i}az del Castillo No. 340,
Colima, Col., Mexico 28045 \email{andres\_pedroza@ucol.mx}}

%\and Name of Second Author \at Name, Address of Institute \email{name@email.address}}
%
% Use the package "url.sty" to avoid
% problems with special characters
% used in your e-mail or web address
%
\maketitle

\abstract{
In this note we present a brief introduction to Lagrangian Floer homology and its relation with the 
solution of Arnol'd conjecture, on the minimal number of non-degenerate fixed points of a Hamiltonian diffeomorphism.
We start with the basic definition of critical point on smooth manifolds, in oder to sketch some aspects of Morse theory.
Introduction to the basics concepts of symplectic geometry are also included, with the idea of understanding the 
statement of Arnol'd Conjecture and how is related to the intersection of Lagrangian submanifolds.
}

%%%%%%%%%%%%%%%%%%%%%%%%%%%%%%%%%%%%%%%%%%%%%%%%%%%%%%%%%%
%%%%%%%%%%%%%%%%%%%%%%%%%%%%%%%%%%%%%%%%%%%%%%%%%%%%%%%%%%
\section{Introduction}
%%%%%%%%%%%%%%%%%%%%%%%%%%%%%%%%%%%%%%%%%%%%%%%%%%%%%%%%%%
%%%%%%%%%%%%%%%%%%%%%%%%%%%%%%%%%%%%%%%%%%%%%%%%%%%%%%%%%%

Many elegant results in mathematics have to deal with the fixed-point-set of a function. For example:
Brouwer fixed-point theorem, Lefschetz fixed-point theorem, Banach fixed-point theorem and
Poincar\'e-Birkhoff theorem, just to name a few. Furthermore, these results 
are fundamental in their own area of mathematics and 
have interesting consequences
in 
diverse areas of mathematics; differential equations, topology and game theory among others. 
 Symplectic geometry has its own fixed-point theorem, which was 
conjectured by V. Arnol'd
\cite{arnold-sur} in 1965. The Arnol'd Conjecture was motivated by
Poincar\'e-Birkhoff theorem: An area-preserving diffeomorphism of the annulus which maps
the boundary circles to themselves in different direction, must have at least two fixed points.

The generalization of Poincar\'e-Birkhoff theorem fits in symplectic geometry and not in
volume-preserving geometry. The Arnol'd Conjecture establishes  a lower  bound on the number
of fixed points a Hamiltonian diffeomorphism in terms of the topology of the manifold.
The fixed points of a Hamiltonian diffeomorphism, (in fact any diffeomorphisms) can be 
seen as the intersection of  its graph and the diagonal. In the context of symplectic geometry, 
 is the intersection of two Lagrangian submanifolds.

In 1987, A. Floer \cite{floer-morse} developed a homological theory that focused on the intersection of Lagrangian
submanifolds. In particular, under some hypotheses, he proved the Arnol'd Conjecture for a particular class
of closed symplectic manifolds.  This theory is called Lagrangian Floer homology.

In these notes we sketch how Lagrangian Floer homology is defined.  In fact we review some aspects of
Morse theory from its basics; like non-degenerate critical points,  the Hessian, flow lines of the
gradient vector field  up to Morse homology. 
The reason being, that Lagrangian Floer homology
emulates in many aspects Morse homology.
Also we cover the basics of symplectic manifolds and Hamiltonian diffeomorphisms. The last section deals
with Lagrangian Floer homology and how it is used to prove the Arnol'd Conjecture.

For the basic notions of differential geometry the reader can look at
\cite{tu-man}; for the aspects of symplectic geometry \cite{cannas-lectures} and \cite{ms}; and also 
\cite{msjholo} where the analytical aspect of holomorphic curves is covered. 
For  details and proofs on the construction  of Lagrangian Floer homology see 
\cite{audin-morseth}, \cite{oh-symplectic1}
and \cite{oh-symplectic2}.
 For an excellent introduction to Fukaya categories see a\cite{auroux-a}; and 
\cite{seidel-fukcat}
for a detail treatment of the subject.

These lecture notes are based on a course
given at the 7th. Mini Meeting on Differential Geometry
held at CIMAT in February 2015. 
The author wishes thank the organizers and participants for
the pleasent atmospere. Finally the author was partially
supported by a CONACYT grant CB-2010/151846.

%%%%%%%%%%%%%%%%%%%%%%%%%%%%%%%%%%%%%%%%%%%%%%%%%%%%%%%%%%
%%%%%%%%%%%%%%%%%%%%%%%%%%%%%%%%%%%%%%%%%%%%%%%%%%%%%%%%%%
\section{Morse-Smale Functions}
%%%%%%%%%%%%%%%%%%%%%%%%%%%%%%%%%%%%%%%%%%%%%%%%%%%%%%%%%%
%%%%%%%%%%%%%%%%%%%%%%%%%%%%%%%%%%%%%%%%%%%%%%%%%%%%%%%%%%

Let $M$ be a smooth manifold of dimension $n$ and $f:M\to \mathbb{R}$  a smooth function.
A point $p\in M$ is called a {\em critical point} of $f$ if the differential $df_p: T_pM\to \mathbb{R}$
at $p$ is the zero map.  Denote by $\textup{Crit}(f)$ the set of critical points of  $f$.
Notice that $\textup{Crit}(f)$ can be the empty set, however  if $M$ is compact then  it  %$\crit{f}$
 is not empty, since a smooth function 
on $M$  has
a maximum and a minimum. 

Let $p\in M$ be a critical point of ${f}$ and  $(x_1,\ldots,x_n)$ a coordinate
chart about $p$. The {\em Hessian matrix} of $f$ at $p$ relative to the chart  $(x_1,\ldots,x_n)$,
is the $n\times n$ matrix 
\begin{eqnarray*}
\textup{Hess}(f,p)= \left (   \frac{\partial^2 f}{\partial x_i \partial x_j}(p)  \right).
\end{eqnarray*}
A critical point $p$ is said to be {\em non-degenerate} if the matrix  $\textup{Hess}(f,p)$
is non-singular. 
Note that the Hessian matrix is symmetric, hence if
it is  non-singular its eigenvalues are real and non-zero. The {\em index}  
of $f$ at a non-degenerate
 critical  point $p$, which is denoted by  $\ind{f}{p}$, 
is defined as the  number of negative eigenvalues of the Hessian matrix at $p$.

The definition of the index at a non-degenerate critical point given above depends
on the coordinate system; however it can be shown that the is independent of the
coordinate system about the the critical point.There is an alternative definition 
of the index of a function at a non-degenerate critical point, that does not needs
 a coordinate system.
%Nonetheless  it is possible to define the index %of a non-degenerate critical point 
%in a coordinate-free manner as follows. 
For  a critical point $p\in M$ of ${f}$ define the bilinear form
$$
df^2_p: T_pM\times T_pM \to\mathbb{R}
$$
as $df^2_p(X,Y):= X(\widetilde Y f)$, where $\widetilde Y$ is
any vector field on $M$
whose value at $p$ is  $Y$. Notice that since $p$ is a critical point of $f$,  the bilinear form
$df^2_p$ is symmetric, 
$$
0=df_p[\widetilde X,\widetilde Y]=[\widetilde X,\widetilde Y]_p(f)=
\widetilde X_p(\widetilde Y f) -\widetilde Y_p(\widetilde X f).
$$
In this context, $p$ is called {\em non-degenerate} if the bilinear symmetric form
$df^2_p$ is non-degenerate. 
%The bilinear form $df^2_p$ induces a symmetric
%non-degenerate  $n\times n$
%matrix.
 The  {\em index} of $f$ at $p$ is defined as the number of negative eigenvalues
of the symmetric bilinear form $df^2_p$.
The two definitions given of non-degenerate critical point 
agree.
The same applies for the two definitions of the  index of a  non-degenerate critical point.
 For further details,
see  \cite[Ch. 1]{audin-morseth} and \cite{liviu-morse}.    %%$\ind{f}{p}$

\begin{definition}
A smooth function $f:M\to \mathbb{R}$ for which all of its critical points
are non-degenerate is called a {\em Morse function.}
\end{definition}

Now we consider some examples in the  case when $M=\mathbb{R}^2$. 
The origin is the only critical point of the function  $f(x,y)=x^2+y^2$.
Moreover is a  non-degenerate critical point and its index is zero. The origin is also the only
 non-degenerate critical point of the functions $g(x,y)=x^2-y^2$ and $h(x,y)=-x^2-y^2$.
In these cases the index at the origin is 1 and 2 respectively.  
These three examples describe
the general behavior of a function on $\mathbb{R}^2$ near the origin when it is a 
non-degenerate critical point.
The precise statement on the behavior of a  function  near a non-degenerate critical point 
is given by Morse lemma. 
\begin{theorem}[Morse Lemma]
Let $f:\mathbb{R}^n\to\mathbb{R}$ be a smooth function such that the origin
is a non-degenerate critical point of index $\lambda$. Then there exists a
coordinate chart $(u_1,\ldots,u_n)$ about the origin such that
$$
f(u_1,\ldots,u_n)= f(0)-u_1^2- \cdots -u_\lambda^2 + u_{\lambda+1}^2 +\cdots +u_n^2.
$$
\end{theorem}

It goes without saying that Morse lemma  also holds  for smooth functions defined
on  arbitrary manifolds. A consequence of Morse lemma, as stated above, is that 
there exists a neighborhood about the origin in $\mathbb{R}^n$ so that it is 
the only critical point in such neighborhood.

\begin{corollary}
Non-degenerate critical points  of a smooth function are isolated.
\end{corollary}

Note that a Morse function defined on a compact manifold has finitely many critical points.

The main reason  behind the study of  Morse functions
is to understand the topology of the manifold. % $M$.
Thus for a smooth function $f:M\to\mathbb{R}$ and 
$a\in \mathbb{R}$ define the {\em level set} 
$$
M_a:=\{x\in M | f(x)\leq a\}\subset M.
$$
Notice that when $a_0$ is the absolute minimum of  $f$, then $M_a$ is empty for every
$a<a_0$. And in the case when $a_1$ is the  absolute maximum of  $f$, then 
$M_a=M$  for every
$a_1\leq a$.

Now we  explain what we mean by  understanding the topology of the manifold;
one aspect is that the manifold can be constructed from  information from
a fixed Morse function on it. 
Consider  a compact manifold  $M$, a smooth Morse function $f:M\to\mathbb{R}$ 
and for simplicity assume that $p_0,\ldots, p_k$ are all the critical points, 
with $\lambda_i=\ind{f}{p_i}$ and
$\lambda_{i}<\lambda_{i+1}$ for $i\in \{0,\ldots, k-1\}$. Thus $f$
achieves its minimum  at $p_0$ and $\ind{f}{p_0}=0$; and 
 it achieves its maximum at   $p_k$ and $\ind{f}{p_k}=n$.
In order to build the manifold $M$ from the critical points of $f$,
one starts 
with the point $M_{f(p_0)}=\{p_0\}$. Then from Theorem \ref{t:morsecell} below, it follows 
that   $M_{a}$ has the same homotopy type has
$M_{f(p_0)}$ for $a\in ( f(p_0),f(p_1) )$.
By an 
{\em $\lambda$-cell} we mean  a space homeomorphic to the closed ball of dimension  $\lambda$.
Hence,   $M_{a}$ is homeomorphic to the $n$-cell for $a\in ( f(p_0),f(p_1) )$.

The next step is to analyze the next non-degenerte critical point $p_1\in M$. 
In this case for $a\in ( f(p_1),f(p_2) )$, it follows that
$M_a$ has the same homotopy  type
as $M_{f(p_0)}$ with a $\lambda_1$-cell  attached. % where $a\in ( f(p_1),f(p_2) )$.
That is $M_a\simeq M_{f(p_0)}\cup_g e_{\lambda_1}$, where 
$g: \partial(e_{\lambda_1}) \to M_{f(p_0)}$ is a gluing function.
This process continues at every critical point.
That is   for $a\in ( f(p_{i}),f(p_{i+1}) )$
the space  $M_a$  as the same homotopy type has
to $ M_{f(p_{i-1})}$ with an attached  $\lambda_{i}$-cell.
The last step asserts that $M_a \simeq  M_{f(p_{k-2})}\cup  e_{\lambda_{k-1}}$
for $a\in ( f(p_{k-1}),f(p_{k}) )$; that is $M_a$  is homeomorphic to $M$ minus an open ball.
Therefore  $M$  is homeomorphic to $M_a$ with a $n$-ball attached. 
Note that the change of topology between the level sets occurs precisely
at the critical points of $f.$
Below, we carry out the same process described above  for $\mathbb{R}P^1.$

Therefore when $f:M\to\mathbb{R}$ is a Morse function, is possible to describe the
topology of the level sets $M_a$ as $a$  increases; in particular  the topology of $M$.
 Furthermore, there is an alternative approach to understand the topology of $M$ using a Morse function.
 This is called Morse homology and it will be describe in Section \ref{s:mh}.

\begin{theorem}
\label{t:morsecell}
Let $f:M\to \mathbb{R}$ be a Morse function.
\begin{itemize}
\item If $f$ has no critical value in $[a,b]$, then $M_a$ is diffeomorphic to $M_b$.
\item If $f$ has only one critical value in  $[a,b]$ of index $\lambda$, then  $M_b$
has the  same  homotopy type as that of $M_a\cup_g  e_\lambda$, for some
gluing function $g$.
\end{itemize}
\end{theorem}
As above, % Here,  $e_\lambda$ stands for a $\lambda$-cell and 
$M_a\cup e_\lambda$ means that
$e_\lambda$ is attached to $M_a$  by some gluing function $g: \partial(e_\lambda) \to M_a$.
Note that $\partial(e_\lambda) $ is diffeomorphic to $S^{\lambda-1}$. In the next example,
we show how Theorem \ref{t:morsecell} is used to obtain the whole manifold $M$, by 
attaching one $\lambda$-cell at a time.

\begin{example}
Consider the real projective space $\mathbb{R}P^n$,  the set of lines through
the origin in $\mathbb{R}^{n+1}$.  A point in $\mathbb{R}P^n$ is represented
in homogeneous coordinates as $[x_0:\ldots :x_n]$.
Let $a_0, \ldots, a_n$ be distinct real numbers, 
define $f: \mathbb{R}P^n\to \mathbb{R}$ by
$$
f([x_0:\cdots :x_n])= \frac{ a_0 x_0^2+ \cdots +a_n x_n^2 }{  x_0^2+ \cdots +x_n^2  }.
$$
So defined
 $f$  is smooth and since the ${a_j}^\prime s$ are distinct it has
$n+1$ non-degenerate critical points, that  are $p_0:=[1:0:\cdots :0], p_1:=[0:1:\cdots :0],\ldots, 
p_n:=[0:\cdots 0:1]$.  Thus $f$ is a Morse function; moreover the critical point $p_j$ has index $j$. 

The reader is encouraged to verify the statements made above. And also
to get the same conclusions for the case of the complex projective space
$\mathbb{C}P^n$ with the function
$$
f([z_0:\cdots :z_n])= \frac{ a_0 |z_0|^2+ \cdots +a_n |z_n|^2 }{  |z_0|^2+ \cdots +|z_n|^2  }.
$$

Now we look at the particular case of $\mathbb{R}P^1$;
recall that   
$\mathbb{R}P^1$ is diffeomorphic to the circle. In this particular case  take 
$a_0=0$ and $a_1=1$, so  $f$ takes the form
$$
f([x_0:x_1])= \frac{   x_1^2 }{  x_0^2 +x_1^2  }.
$$
In this case
$f$ has only one critical point 
of index $0$, namely at $[1:0]$. It also has only one critical point  of index $1$, at $[0:1]$. 
These points correspond  to the maximum and minimum of $f$.
In terms of Theorem	\ref{t:morsecell} the circle is obtained as follows.  We start with the
$0$-cell that is just a point, that is $M_0=\{[1:0]\}$. Since $f$ has no critical values
in the interval $[0,1/2]$ other than $0$, then from Theorem \ref{t:morsecell} if follows
that $M_0$ has the same homotopy type has $M_{1/2}$. Notice that $M_{1/2}$ is a semicircle,
the south hemisphere. Next comes the other critical point $[0:1]$.
 It has index $1$; thus a $1$-cell
is attached to $M_{1/2}$. That is, the two points of $\partial (e_1)$ get glued to     
$M_{1/2}$ to obtain the circle.
See Figure \ref{f:morsecircle}.
\begin{figurehere}
\begin{center}
{\psset{unit=.8}
\begin{pspicture}(-6,-3)(8,3) %%%% Mover el 10 para mover la figura 
\pscircle(6,0){2}
%\rput(0,2.5){\parametricplot{135}{415}{t cos .5 mul t sin .15 mul}}
%\parametricplot[linestyle=dashed]{0}{80}{t cos 2 mul t sin .6 mul}
%\rput(0,2.5)
\parametricplot{-1.5}{1.5}{t 2}
\parametricplot{180}{360}{t cos 2 mul t sin 2 mul}
\psdots(-5,-2)(-1.5,2)(1.5,2) 
\psdots(-2,0)(-1.5,2)(1.5,2)  %(0,-2) (0,2)
\psdots(2,0)(-1.5,2)(1.5,2)  %(0,-2) (0,2)
\uput{.2}[0](-5.8,-1.3){$0$-cell}
\uput{.2}[0](-.5,-1.3){$M_{1/2}$}
\uput{.2}[0](-.6,1.4){$1$-cell}
\uput{.2}[0](7.3,1.7){$\mathbb{R}P^1$}
\uput{.2}[0](-5.4,-2.7){(a)}
\uput{.2}[0](-.4,-2.7){(b)}
\uput{.2}[0](5.6,-2.7){(c)}
\end{pspicture}
}
\caption{ Morse decomposition of $\mathbb{R}P^1$ with respecto to $f$.
 In (a) the $0$-cell that corresponds to $p_0$
of index $0$. In (b) the submanifold  $M_{1/2}$, diffeomorphic to a point, is attached a $1$-cell that
corresponds to the point $p_1$ of index 1.
Finally,  (c) the result after attaching the $1$-cell. }
\label{f:morsecircle}
\end{center}
\end{figurehere}

\end{example}
\medskip
An important aspect to consider is the existence of Morse functions on a given manifold. 
It turns out that there are plenty of Morse functions. More precisely, the set  
of Morse functions on a closed manifold is $C^2$-dense in the space of smooth functions. The reason that
the $C^2$-topology is needed is because  the concept of non-degenerate critical points
involves  derivatives up to second-order. 
In theory, is not to difficult to understand the topology of $M$ via a 
Morse function as above. Next we take this idea a step further to recover the homology of $M$.

Fix a Riemannian metric $g$  on $M$ and let $\langle\cdot,\cdot\rangle $ be the induced inner product
on its tangent bundle. 
The {\em gradient
vector field}, $\textup{grad}(f)$, of the function $f:M\to \mathbb{R}$ is defined
 by the equation
$$
\langle \grad{f}, X\rangle =X(f)
$$
for every vector field $X$ on $M$. Notice that if $p$ is a critical point of $f$, then
$\grad{f}_p=0$. And  conversely, if  $\grad{f}_p=0$  then
$p$ is a critical point of $f.$  Therefore $\crit{f}$ equals  the zero set of
$\grad{f}$.

In order to simplify the exposition, from now on we assume that $M$ is compact.
Denote by $\theta:\mathbb{R}\times M\to M$ the flow of the {\bf negative} 
gradient vector field of $f$. Thus for $x\in M$
$$
\left.\frac{\partial\theta(t,x)}{\partial t}\right|_{t=0}=-\grad{f}_x.
$$
The reason to consider the negative gradient vector field is only a matter of convention.
Note that $-\grad{f}=\grad{-f}$ and if $p$ is a non-degenerate critical point of $f$,
then $\ind{f}{p} = n-\ind{-f}{p}$ where
$n$ is the dimension of $M$.
Also notice that $-\grad{f}(f)<0$ outside the set of critical points of
$f$, hence $-\grad{f}$ points in the direction in which $f$ is decreasing.
The way to think about the index of a non-degenerate critical point is the
number of linearly independent directions in which the $-\grad{f}$ decreases.
Let $p$ be a point where $-\grad{f}$ vanishes, then 
consider  all points of $M$ that under the flow $\theta$ converge to
$p$ as $t$ goes to infinity;
$$
W^s(f,p):=\left\{x\in M \left |   \lim_{t\to+\infty} \theta(t,x)=p \right.\right\} .
$$
Similarly, % we have the unstable submanifold at $p$
$$
W^u(f,p):=\left\{x\in M \left |  \lim_{t\to-\infty} \theta(t,x)=p \right. \right\},
$$
the set of all points in $M$ that have $p$ has a source. 
Since $-\grad{f}$ vanishes at $p$, then the critical point $p$
is a fixed under the flow, hence $p\in W^s(f,p)$ and $p\in W^u(f,p)$.
The submanifolds $W^u(f,p)$ and $W^s(f,p)$ are called the {\em unstable manifold} and 
{\em stable
submanifold} of $f$ at $p$, respectively.

%Since $p$ is a fix point of the flow, then $(\theta^{t})_{*,p}:T_pM\to T_pM $ 
%for every $t$. Note that $p$ is a hyperbolic point of the flow of $-\grad{f}$
%if and  only if $p$ is a non-degenerate point of $f$

\begin{theorem}
\label{t:unstable}
If $p$ is a non-degenerate critical point of $f$, then $W^u(f; p)$
is a smooth submanifold of $M$ of dimension $\ind{f}{p}$.
\end{theorem}

Instead, if we consider the function $-f$  the set critical non-degenerate
points of $f$ and $-f$ agree. Moreover $W^u(f; p)=W^s(-f; p)$ and 
$W^s(f; p)=W^u(-f; p)$. Hence $W^s(f; p)$ is also a smooth submanifold
of $M$ of dimension $n-\ind{f}{p}$. 

\begin{example}
\label{exa:sphere}
Let $M=S^2$ be the unit sphere in $\mathbb{R}^3$ centered at the origin and
$f:S^2\to \mathbb{R}$ defined as $f(x,y,z)=z$. Then the poles
$N=(0,0,1)$ and $S=(0,0,-1)$ are  the critical points of $f$. Furthermore they are
non-degenerate,
$N$ has index 2 and $S$ has index 0. 

\begin{figurehere}
\begin{center}
{\psset{unit=.8}
\begin{pspicture}(-2,-2.6)(2,3)
%\pscircle(0,0){2}
\parametricplot[linewidth=.5pt]{0}{180}{t sin .3 mul t cos 2 mul} 
\parametricplot[linewidth=.5pt]{0}{180}{t sin .8 mul t cos 2 mul} 
\parametricplot[linewidth=.5pt]{0}{180}{t sin 1.3 mul t cos 2 mul} 
\parametricplot[linewidth=.5pt]{0}{180}{t sin 1.8 mul t cos 2 mul} 

\parametricplot[linewidth=.5pt]{0}{180}{t sin -.3 mul t cos 2 mul} 
\parametricplot[linewidth=.5pt]{0}{180}{t sin -.8 mul t cos 2 mul} 
\parametricplot[linewidth=.5pt]{0}{180}{t sin -1.3 mul t cos 2 mul} 
\parametricplot[linewidth=.5pt]{0}{180}{t sin -1.8 mul t cos 2 mul} 

\parametricplot[linewidth=.5pt,arrows=->]{0}{40}{t sin .3 mul t cos 2 mul}
\parametricplot[linewidth=.5pt,arrows=->]{40}{60}{t sin .3 mul t cos 2 mul}
\parametricplot[linewidth=.5pt,arrows=->]{60}{80}{t sin .3 mul t cos 2 mul}
\parametricplot[linewidth=.5pt,arrows=->]{80}{100}{t sin .3 mul t cos 2 mul}
\parametricplot[linewidth=.5pt,arrows=->]{100}{120}{t sin .3 mul t cos 2 mul}
\parametricplot[linewidth=.5pt,arrows=->]{120}{140}{t sin .3 mul t cos 2 mul}

\parametricplot[linewidth=.5pt,arrows=->]{0}{40}{t sin .8 mul t cos 2 mul}
\parametricplot[linewidth=.5pt,arrows=->]{40}{60}{t sin .8 mul t cos 2 mul}
\parametricplot[linewidth=.5pt,arrows=->]{60}{80}{t sin .8 mul t cos 2 mul}
\parametricplot[linewidth=.5pt,arrows=->]{80}{100}{t sin .8 mul t cos 2 mul}
\parametricplot[linewidth=.5pt,arrows=->]{100}{120}{t sin .8 mul t cos 2 mul}
\parametricplot[linewidth=.5pt,arrows=->]{120}{140}{t sin .8 mul t cos 2 mul}

\parametricplot[linewidth=.5pt,arrows=->]{0}{40}{t sin 1.3 mul t cos 2 mul}
\parametricplot[linewidth=.5pt,arrows=->]{40}{60}{t sin 1.3 mul t cos 2 mul}
\parametricplot[linewidth=.5pt,arrows=->]{60}{80}{t sin 1.3 mul t cos 2 mul}
\parametricplot[linewidth=.5pt,arrows=->]{80}{100}{t sin 1.3 mul t cos 2 mul}
\parametricplot[linewidth=.5pt,arrows=->]{100}{120}{t sin 1.3 mul t cos 2 mul}
\parametricplot[linewidth=.5pt,arrows=->]{120}{140}{t sin 1.3 mul t cos 2 mul}

\parametricplot[linewidth=.5pt,arrows=->]{0}{40}{t sin 1.8 mul t cos 2 mul}
\parametricplot[linewidth=.5pt,arrows=->]{40}{60}{t sin 1.8 mul t cos 2 mul}
\parametricplot[linewidth=.5pt,arrows=->]{60}{80}{t sin 1.8 mul t cos 2 mul}
\parametricplot[linewidth=.5pt,arrows=->]{80}{100}{t sin 1.8 mul t cos 2 mul}
\parametricplot[linewidth=.5pt,arrows=->]{100}{120}{t sin 1.8 mul t cos 2 mul}
\parametricplot[linewidth=.5pt,arrows=->]{120}{140}{t sin 1.8 mul t cos 2 mul}

\parametricplot[linewidth=.5pt,arrows=->]{0}{40}{t sin 2 mul t cos 2 mul}
\parametricplot[linewidth=.5pt,arrows=->]{40}{60}{t sin 2 mul t cos 2 mul}
\parametricplot[linewidth=.5pt,arrows=->]{60}{80}{t sin 2 mul t cos 2 mul}
\parametricplot[linewidth=.5pt,arrows=->]{80}{100}{t sin 2 mul t cos 2 mul}
\parametricplot[linewidth=.5pt,arrows=->]{100}{120}{t sin 2 mul t cos 2 mul}
\parametricplot[linewidth=.5pt,arrows=->]{120}{140}{t sin 2 mul t cos 2 mul}
\parametricplot[linewidth=.5pt,arrows=->]{140}{180}{t sin 2 mul t cos 2 mul}

\parametricplot[linewidth=.5pt,arrows=->]{0}{40}{t sin -.3 mul t cos 2 mul}
\parametricplot[linewidth=.5pt,arrows=->]{40}{60}{t sin -.3 mul t cos 2 mul}
\parametricplot[linewidth=.5pt,arrows=->]{60}{80}{t sin -.3 mul t cos 2 mul}
\parametricplot[linewidth=.5pt,arrows=->]{80}{100}{t sin -.3 mul t cos 2 mul}
\parametricplot[linewidth=.5pt,arrows=->]{100}{120}{t sin -.3 mul t cos 2 mul}
\parametricplot[linewidth=.5pt,arrows=->]{120}{140}{t sin -.3 mul t cos 2 mul}

\parametricplot[linewidth=.5pt,arrows=->]{0}{40}{t sin -.8 mul t cos 2 mul}
\parametricplot[linewidth=.5pt,arrows=->]{40}{60}{t sin -.8 mul t cos 2 mul}
\parametricplot[linewidth=.5pt,arrows=->]{60}{80}{t sin -.8 mul t cos 2 mul}
\parametricplot[linewidth=.5pt,arrows=->]{80}{100}{t sin -.8 mul t cos 2 mul}
\parametricplot[linewidth=.5pt,arrows=->]{100}{120}{t sin -.8 mul t cos 2 mul}
\parametricplot[linewidth=.5pt,arrows=->]{120}{140}{t sin -.8 mul t cos 2 mul}

\parametricplot[linewidth=.5pt,arrows=->]{0}{40}{t sin -1.3 mul t cos 2 mul}
\parametricplot[linewidth=.5pt,arrows=->]{40}{60}{t sin -1.3 mul t cos 2 mul}
\parametricplot[linewidth=.5pt,arrows=->]{60}{80}{t sin -1.3 mul t cos 2 mul}
\parametricplot[linewidth=.5pt,arrows=->]{80}{100}{t sin -1.3 mul t cos 2 mul}
\parametricplot[linewidth=.5pt,arrows=->]{100}{120}{t sin -1.3 mul t cos 2 mul}
\parametricplot[linewidth=.5pt,arrows=->]{120}{140}{t sin 1.3 mul t cos 2 mul}

\parametricplot[linewidth=.5pt,arrows=->]{0}{40}{t sin -1.8 mul t cos 2 mul}
\parametricplot[linewidth=.5pt,arrows=->]{40}{60}{t sin -1.8 mul t cos 2 mul}
\parametricplot[linewidth=.5pt,arrows=->]{60}{80}{t sin -1.8 mul t cos 2 mul}
\parametricplot[linewidth=.5pt,arrows=->]{80}{100}{t sin -1.8 mul t cos 2 mul}
\parametricplot[linewidth=.5pt,arrows=->]{100}{120}{t sin -1.8 mul t cos 2 mul}
\parametricplot[linewidth=.5pt,arrows=->]{120}{140}{t sin -1.8 mul t cos 2 mul}

\parametricplot[linewidth=.5pt,arrows=->]{0}{40}{t sin -2 mul t cos 2 mul}
\parametricplot[linewidth=.5pt,arrows=->]{40}{60}{t sin -2 mul t cos 2 mul}
\parametricplot[linewidth=.5pt,arrows=->]{60}{80}{t sin -2 mul t cos 2 mul}
\parametricplot[linewidth=.5pt,arrows=->]{80}{100}{t sin -2 mul t cos 2 mul}
\parametricplot[linewidth=.5pt,arrows=->]{100}{120}{t sin -2 mul t cos 2 mul}
\parametricplot[linewidth=.5pt,arrows=->]{120}{140}{t sin -2 mul t cos 2 mul}
\parametricplot[linewidth=.5pt,arrows=->]{140}{180}{t sin -2 mul t cos 2 mul}
\uput{.2}[0](-.3,-2.3){S}
\uput{.2}[0](-.3,2.3){N}
\end{pspicture}
}
\caption{The flow lines of the gradient vector field  of $f(x,y,z)=z$ on $S^2$ with respect to the
standard Riemannian structure.}
\label{f:hopf}
\end{center}
\end{figurehere}

Consider  the Riemannian structure on $S^2$ induced from the standard Riemannian structure
on $\mathbb{R}^3$. Then $-\textup{grad}(f)$ is the vector field that points downwards, 
and 
\begin{eqnarray*}
\label{eq:s2us}
W^u(f,N)=S^2\setminus \{S\}, W^s(f,N)=\{N\}, W^u(f,S)=\{S\}
\mbox{ and } W^s(f,S)= S^2\setminus\{N\}.
\end{eqnarray*}
\end{example}

\begin{example}
\label{exa:torus}
Consider the function $f:\mathbb{R}^2\to \mathbb{R}$ given by
$$
f(x,y)=\cos(2\pi x)+ \cos(2\pi y).
$$
So defined $f$ induces a smooth function on the flat two-dimensional torus 
$\mathbb{T}^2=\mathbb{R}^2/\mathbb{Z}^2$,
which we still denote by $f$. There are $4$ non-degenerate critical points on the torus,
$p_1=[0,0], p_2=[0,1/2],p_3=[1/2,0]$ and $p_4=[1/2,1/2],$   of index
$2,1,1$ and $0$ respectively.  Consider the Riemannian structure on $\mathbb{T}^2$
induced from the canonical Riemannian structure on $\mathbb{R}^2$. 
Then the  flow of $-\textup{grad}(f)$ can be seen in  Figure
\ref{f:flowtorus}.

\begin{figurehere}
\begin{center}
\includegraphics[scale=0.5]{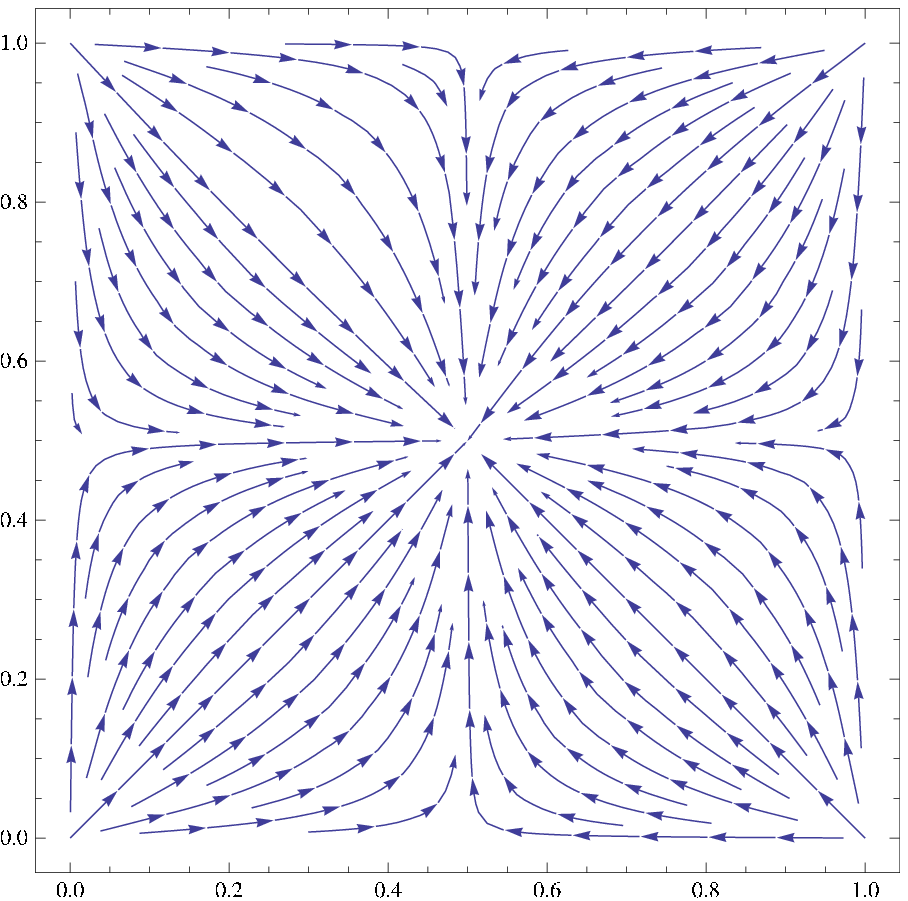}
\caption{
The flow lines of the gradient vector field  of $\cos(2\pi x)+ \cos(2\pi y)$ on  $\mathbb{T}^2$ with respect 
to the flat Riemannian structure.}
\label{f:flowtorus}
\end{center}
\end{figurehere}

Notice that there are only  two lines that connect $p_2$ to $p_4$. And a 1-dimensional family of 
flow lines that connect $p_1$ to $p_4$, whose points determine four open connected 
components of the torus.

Also  Figure
\ref{f:flowtorus} gives a description of the stable and unstable submanifolds.
Observe  that every interior point of $[0,1]\times [0,1]$, lies in a flow
line that ends at $p_4$. That is, 
$$
W^s(f; p_4)=\mathbb{T}^2 \setminus \{ \partial([0,1]\times [0,1])\}.
$$
Similarly we have that $W^u(f; p_1)$ equals
$$
 \mathbb{T}^2 \setminus \{  p_2,p_3,p_4  \} \cup
\{ (x,1/2) |  x\in [0,1]\setminus \{1/2\}\}\cup \{ (1/2,y) |  y\in [0,1]\setminus \{1/2\}\},
$$
and
\begin{eqnarray*}
W^u(f; p_2) &=&   \{ (x,1/2) |  x\in [0,1]\setminus \{1/2\}\},\\
W^u(f; p_3)&=&   \{ (1/2,y) |  y\in [0,1]\setminus \{1/2\}\}.
\end{eqnarray*}
\end{example}

%\begin{figure}[h!]
%\centering
%\includegraphics[scale=0.5]{torus}
%\centering
%
%
%\end{figure}

Let $p$ and $q$ be non-degenerate critical points of a smooth function $f:M\to\mathbb{R}$.
Then  the set $W^u(f; p)\cap W^s(f; q) $ consists of points of $M$ that belong to a  flow line
$u:\mathbb{R}\to M$ of $-\grad{f}_{u(t)}$ that connects $p$ to $q$; that is
\begin{eqnarray}
\label{e:flowlines}
\frac{d u}{dt}(t)=-\grad{f}_{u(t)}, \ \
\lim_{t\to-\infty} u(t)=p  \ \ \textup{ and }\ \
\lim_{t\to+\infty} u(t)=q. 
\end{eqnarray}
We know from Theorem \ref{t:unstable} that $W^u(f; p)$  and $W^s(f; q)$ are submanifolds
of $M$, but their intersection might not be a smooth manifold. Hence
a  smooth function $f:M\to\mathbb{R}$ is said to satisfy the {\em Smale condition}
if for any pair of critical points
$p$ and $q$, $W^u(f; p)$ and $W^s(f; q) $ intersect transversally. 
In particular $W^u(f; p)\cap W^s(f; q) $ is a submanifold of $M$.

The function that appears in Example \ref{exa:sphere} satisfies the Smale condition.
In this example the intersection of any  pair of stable and unstable submanifolds is either
empty, a point, the sphere minus a point or  the sphere minus two points.
Also the Moorse function in Example \ref{exa:torus} satisfies the Smale condition.
In particular, notice that $W^u(f; p_2)\cap W^s(f; p_4)$ consists of
two disjoint open intervals.

The type of functions that are of interest in this note are the Morse-Smale functions.
For an arbitrary compact manifold and Riemannian metric, there always exists a Morse-Smale
function. Furthermore, in some sense there are plenty of such functions. Then if $(M, g)$
is a Riemannian manifold and $f:M\to\mathbb{R}$  a Morse-Smale function
 we write
$\mathcal{M}(f;p,q)$ for the set of points of $M$ that belong to a  flow trajectory
of $-\grad{f}$
that goes  from $p$ to $q$ as in Eq. (\ref{e:flowlines}).
Notice that in this case 
$\mathcal{M}(f;p,q)$ is a smooth submanifold of $M$ of dimension
$\ind{f}{p}-\ind{f}{q}$.
Note that the submanifold $\mathcal{M}(f;p,q)$ admits a natural   action of $\mathbb{R}$ defined
as $(s.u)(t):=u(t+s)$ for $s\in \mathbb{R}$. The action is in fact free and the orbit space
of this action is denoted by $\hat{\mathcal{M}}(f;p,q)$. Hence $\hat{\mathcal{M}}(f;p,q)$ is
identified as the space of trajectories that joint $p$ to $q.$

So defined,  the space of points that belong to a flow line of $-\grad{f}$ that connect
$p$ to $q$, $\mathcal{M}(f;p,q)$, is not necessarily compact. For example, in the case of the 
two-spere in Example \ref{exa:sphere} we have that $\mathcal{M}(f;N,S)$ is $S^2\setminus\{N,S\}$.
In this example if we add the critical points we obtain a compact space, namely the whole manifold
$S^2$. Note that in this example $\hat{\mathcal{M}}(f;N,S)$ is diffeomorphic to $S^1.$
But it is not always the case that by adding the critical points $p$ and $q$ to
$\mathcal{M}(f;p,q)$ that it becomes a compact space.
For instance, in the torus case of Example \ref{exa:torus} the space $\mathcal{M}(f;p_1,p_4)\cup \{p_1,p_4\}$
is not compact.

In general, the way to compactify the space of trajectories $\hat{\mathcal{M}}(f;p,q)$
 is by adding broken trajectories. A {\em
broken trajectory} from $p$ to $q$ is a collection of flow lines $\{u_1,\ldots, u_r \}$
of $-\grad{f}$ such that $u_j$ connects the critical points $x_{j}$ to $x_{j+1}$ for 
$j\in\{1,\ldots ,r\}$ 
where $p=x_1$ and   $q=x_{r+1}.$  Consider the bigger set of flow lines that
connect $p$ to $q$, namely usal flow trajectories plus   broken trajectories,
$$
\overline{\mathcal{M}}(f;p,q):=\hat{\mathcal{M}}(f;p,q)\cup\{\mbox{broken 
trajectories from $p$ to $q$} \}.
$$
Recall that the index of critical points of $f$ decreases along flow lines. Hence
the number of flow lines that form a broken flow lines is less than 
$\ind{f}{p}-\ind{f}{q}$. Hence if $\ind{f}{p}-\ind{f}{q}=1$ there are no
broken trajectories connecting $p$ to $q$ and 
$\overline{\mathcal{M}}(f;p,q)=\hat{\mathcal{M}}(f;p,q)$. That is, $\hat{\mathcal{M}}(f;p,q)$
is compact in this case and it consists of finitely many points.

The proof of the  next  result  can consulted in   \cite[Chp. 3]{audin-morseth} and 
\cite{salamon-lectures}.

\begin{proposition}
\label{p:moduli}
Let $(M, g)$ be a closed Riemannian manifold, $f:M\to\mathbb{R}$ a Morse-Smale function and $p,q$
critical points of $f$.
Then the   natural action of $\mathbb{R}$   on $\mathcal{M}(f;p,q)$ is free. Moreover 
$\overline{\mathcal{M}}(f;p,q)$ %:=\mathcal{M}(f;p,q)/\mathbb{R}$
is smooth and compact of dimension $\ind{f}{p}-\ind{f}{q}-1.$
\end{proposition}

The important case that would be relevant later on is the case when $\ind{f}{p}-\ind{f}{q}=2$.
Usually in this case the space $\hat{\mathcal{M}}(f;p,q)$ is not compact, so we must add 
broken trajectories. Hence  $\overline{\mathcal{M}}(f;p,q)$ is a finite collection of closed intervals
and circles.

In Example \ref{exa:torus} consider $u_1(t):=[0, t]$ and $u_2(t):=[t,1/2]$, for $t\in (0, 1/2),$ two
flow lines of $-\grad{f}$. The flow line $u_1$ connects $p_1$ to $p_2$, and 
 $u_2$ connects $p_2$ to $p_4$. Hence $\{u_1,u_2\}$ is a broken trajectory that connects
$p_1$ to $p_4$. Note that $\hat{\mathcal{M}}(f;p_1,p_4)$ is diffeomorphic to
four copies of $(0,1)$; and there are eight broken trajectories that must be added
to obtain $\overline{\mathcal{M}}(f;p,q)$. For instance, $\{u_1,u_2\}$
is one of them. Henceforth $\overline{\mathcal{M}}(f;p,q)$ is diffeomorphic to 
four copies of $[0,1].$

%%%%%%%%%%%%%%%%%%%%%%%%%%%%%%%%%%%%%%%%%%%%%%
%%%%%%%%%%%%%%%%%%%%%%%%%%%%%%%%%%%%%%%%%%%%%% 
\section{Morse Homology}
\label{s:mh}
%%%%%%%%%%%%%%%%%%%%%%%%%%%%%%%%%%%%%%%%%%%%%%
%%%%%%%%%%%%%%%%%%%%%%%%%%%%%%%%%%%%%%%%%%%%%%

We are  going to define the Morse-Witten complex of $(M;f,g)$; the Riemannian manifold 
and the Morse-Smale function. 
For simplicity we will
use $\mathbb{Z}_2$ coefficients, keep in mind that it is possible to use 
integer coefficients.
In order to define Morse homology with integer coefficients, one must prove that
is possible to have a coherent system of orientation on the compact moduli spaces.
In the case of $\mathbb{Z}_2$ coefficients the orientation of the moduli spaces
is irrelevant, only the boundary components of the moduli spaces of dimension two
are important.
See for example \cite{salamon-lectures}, where they use integer coefficients.
Also we drop the dependence of the Riemannian metric from the notation.
Denote by 
 $\textup{Crit}_\lambda(f) $ the set of critical points of index $\lambda$
and by $C_\lambda(f)$ the $\mathbb{Z}_2$-vector space generated by the elements of
 $\textup{Crit}_\lambda(f) $.
For $\lambda\notin \{0,1,\ldots n\}$, define $C_\lambda(f)$ to be the trivial vector space.
If $p$ and $q$ are critical points of $f$ such that $\ind{f}{p}=\ind{f}{q}+1$, then by Proposition
\ref{p:moduli} $\hat{\mathcal{M}}(f;p,q)$ is a finite set of points.
Denote by 
$\#_{\mathbb{Z}_2} \hat{\mathcal{M}}(f;p,q) $ the number of points of
$\hat{\mathcal{M}}(f;p,q)$ module 2.

The boundary operator, $\partial_\lambda: C_\lambda(f)\to C_{\lambda-1}(f)$, is the linear map
 defined on generators $p\in  C_\lambda(f)$ as
\begin{eqnarray*}
%\label{e:partial}
\partial_\lambda (p):=\sum_{q\in \textup{Crit}_{  \lambda-1}(f)}  
\#_{\mathbb{Z}_2} \hat{\mathcal{M}}(f;p,q) \, q.
\end{eqnarray*}
%where $p$ is a critical point of index $\lambda$. 
Notice that if $\hat{\mathcal{M}}(f;p,q)$
is zero-dimensional, then ${\mathcal{M}}(f;p,q)$ consists of finitely many lines
that connect $p$ to $q$. This geometric description of ${\mathcal{M}}(f;p,q)$ is 
 useful when computing the boundary operator
$\partial$; this will be seen  for instance  below in Example \ref{e:homologytorus}.
The reason why $ \partial_{\lambda}$ is called the boundary operator is given by the next result.

In order to compute $ \partial_{\lambda-1}\circ \partial_\lambda$ one must consider
the moduli spaces $\overline{\mathcal{M}}(f;p,r)$ 
where  $\ind{f}{p}-\ind{f}{r}=2.$  For $p \in \textup{Crit}_{  \lambda-2}(f)$,
$$
\partial_{\lambda-1}\partial_\lambda (p):=
\sum_{r\in \textup{Crit}_{  \lambda-2}(f)}  
\sum_{q\in \textup{Crit}_{  \lambda-1}(f)}  
\#_{\mathbb{Z}_2} ( \#\hat{\mathcal{M}}(f;p,q) \cdot \, 
\# \hat{\mathcal{M}}(f;q,r) )\, r
$$
where $\# \hat{\mathcal{M}}(f;p,q)$ stands for the number of points of $\hat{\mathcal{M}}(f;p,q)$. 
Notice that 
$\overline{\mathcal{M}}(f;p,r)$  is a one-dimensional compact manifold; hence it is the union of a finite
collection of closed intervals and circles. Hence its boundary consists of a even number of points
which are
$$
\cup_{q\in \textup{Crit}_{  \lambda-1}(f)}  \#\hat{\mathcal{M}}(f;p,q) \cdot \, 
\# \hat{\mathcal{M}}(f;q,r) 
$$ 
and correspond to the broken trajectories from $p$ to $r$ that go thru $q$.

\begin{theorem}
The  operator satisfies $ \partial_{\lambda-1}\circ \partial_\lambda =0$.
\end{theorem}

The complex $(%C_*(f)=
\oplus_\lambda C_\lambda(f), \partial)$ is called the {\em Morse-Witten
complex} of $(M;f,g)$. Its homology
$$
\textup{MH}_\lambda(M; f,g): =\frac{\textup{Ker}\,  \partial_{\lambda}}{\textup{Im}\,\partial_{\lambda+1} }
$$
is called the Morse homology of $(M; f,g)$ with $\mathbb{Z}_2$-coefficients. 
Note that the relevant moduli spaces $\overline{\mathcal{M}}(f;p,q)$ for the definition of 
Morse homology are those whose dimension is at most two.

%\begin{comment}
\begin{remark}
As mention above is possible to define Morse homology with $\mathbb{Z}$ coefficients.
For, $W^u(f,p)$ is an  orientable submanifold of $M$ for any critical point $p$. Hence,
one fixes an orientation on $W^u(f,p)$ for every critical point. This yields an orientation
on $W^u(f,p)\cap W^s(f,q)$ and hence on ${\mathcal{M}}(f;p,q) $ and $\hat{\mathcal{M}}(f;p,q)$.
Thus if $\ind{f}{p}=\ind{f}{q}+1$, then $\overline{\mathcal{M}}(f;p,q) =
\hat{\mathcal{M}}(f;p,q) $ is a finite set
of points each of which has   a sign. Set $n(f;p,q)$ to be the sum of these signs;  then the boundary operator
over $\mathbb{Z}$-coefficients is defined as
$$
\partial_\lambda (p):=\sum_{q\in \textup{Crit}_{  \lambda-1}(f)}  
n(f;p,q) \, q.
$$
However the statement $\partial_{\lambda-1} \circ \partial_\lambda =0$ is delicate in this 
case. 
%However in order to guarantee that 
One must take into consideration that 
the orientation of the moduli spaces of dimension two induced the right 
orientation on its boundary; the  one-dimensional  moduli spaces. For example see \cite{salamon-lectures}. 
\end{remark}
%\end{comment}

Recall from Example \ref{exa:sphere}, that on $S^2$ we defined
a Morse function with the poles $N$ and $S$ as critical points of index 2 and 0 respectively.
The Riemannian structure on the sphere was induced from the canonical Riemannian structure on 
$\mathbb{R}^3$. 
Further, we calculated the stable and unstable submanifolds of $N$ and $S$.
%Eq. (\ref{eq:s2us}).
%W^u(f,N)=S^2\setminus \{S\} \hbox{ and }W^s(f,N)=\{N\}.%
%$$
From this calculation, it follows that $f$ is a Morse-Smale function. 
Therefore $C_0(f)=\mathbb{Z}_2 \langle S\rangle$, $C_2(f)=\mathbb{Z}_2\langle N\rangle$
and the boundary operator is the zero map. Hence
%$MH_\lambda(S^2; f,g)$ equals  $\mathbb{Z}_2 $ for $\lambda=0,2$ and is trivial otherwise.
$$
\textup{MH}_\lambda(S^2; f,g) = 
\begin{cases} 
\mathbb{Z}_2 & \quad \text{if } \lambda=0,2\\ 
0 &\quad \text{if } \lambda\neq 0,2.\\ 
 \end{cases}
$$
\begin{example}
\label{e:homologytorus}
In this example we consider the function on the two-dimensional torus
$\mathbb{T}^2$ defined in Example \ref{exa:torus}.  Notice that the function is Morse-Smale. Hence
$C_0(f)=\mathbb{Z}_2\langle p_4\rangle, C_1(f)=\mathbb{Z}_2\langle p_2, p_3\rangle,$ and
$C_2(f)=\mathbb{Z}_2\langle p_1\rangle$.
Counting trajectory flow lines, we get
$\partial p_1= 2p_2+2p_3=0$, $\partial p_2=0$, $\partial p_3=0$, and $\partial p_4=0$.
Therefore,
$$
\textup{MH}_\lambda(\mathbb{T}^2; f,g) = 
\begin{cases} 
\mathbb{Z}_2 & \quad \text{if } \lambda=0,2\\ 
\mathbb{Z}_2\times \mathbb{Z}_2 & \quad \text{if } \lambda=1\\
0 &\quad \text{if } \lambda\neq 0,1,2\\ 
 \end{cases}
$$
\end{example}

Summing up, we started with a smooth closed manifold $M$, then we
choose a smooth function $f$ and a Riemannian metric $g$, such that the critical points
of $f$ were non-degenerate and the intersection of the stable and unstable submanifolds
were transversal. With all these data, we defined the Morse homology of $(M;f,g)$ 
At the end, Morse homology is a topological invariant of the manifold; that is, 
is independent of the function and the Riemannian metric. Furthermore 
it recovers the ordinary homology of the manifold.  See \cite{schwarz-mor} and \cite{witten-super}.

\begin{theorem}
\label{t:morse}
Let $M$ be a compact manifold, $g$ a Riemannian metric and $f$ a Morse-Smale
function. Then $
\textup{MH}_*(M; f,g)$ is independent of the function and the metric. Moreover
$\textup{MH}_\lambda(M; f,g)\simeq 
H_\lambda(M; \mathbb{Z}_2)$ as vectors spaces for every $\lambda$.
\end{theorem}

In the words of R. Bott \cite{bott-morse}, {\em Morse theory indomitable}. 
Here we barely treated the subject and its consequences.
The reader is encouraged to learn more about the subject in \cite{bott-morse}, 
\cite{matsumoto-morse}, \cite{liviu-morse} and in the beautiful monograph of J. Milnor, \cite{milnor-morse}.
One application of Morse theory is the handlebody decomposition of a manifold; a much finer 
result than that stated in  Theorem \ref{t:morse}. In particular the Bott periodicity theorem is a marvelous 
consequence of Morse theory.
Another typical consequence of Morse theory are the Morse inequalities. Here the problem is to determine
lower bounds for the number of critical points of a fixed index of a Morse function.
Denote by $b_\lambda(M)$ the 
$\lambda$-Betti number of $M$, that is the rank of $ H_\lambda(M;\mathbb{Z}_2)$.

\begin{theorem}[Morse's inequalities]
\label{t:morseine}
Let $M$ be a a closed manifold and $f:M\to \mathbb{R}$  a Morse
function. Then 
$$
\#\textup{ Crit}_\lambda(f)\geq  b_\lambda(M)
$$
for every  $0\leq \lambda\leq n$.
\end{theorem}

Thus Morse theory gives  a lower bound for the minimal number of critical
points that a Morse function can have on a manifold. 
Finally we mention that some features of ordinary homology, for example Poincar\'e
duality and product operations, can be described in the Morse homology setting.   
See \cite{fukaya-morsehomotopy},
\cite{schwarz-mor} and \cite{witten-super}.

%%%%%%%%%%%%%%%%%%%%%%%%%%%%%%%%%%%%%%%%%%%%%%
%%%%%%%%%%%%%%%%%%%%%%%%%%%%%%%%%%%%%%%%%%%%%%
\section{Symplectic Manifolds and Lagrangian Submanifolds}
\label{s:symandlag}
%%%%%%%%%%%%%%%%%%%%%%%%%%%%%%%%%%%%%%%%%%%%%%
%%%%%%%%%%%%%%%%%%%%%%%%%%%%%%%%%%%%%%%%%%%%%%

A {\em symplectic form} on a manifold $M$ is $2$-form  $\omega$  that is  closed $d\omega=0$ and  non-degenerate.
Here non-degenerate means that at every $p\in M$ and every nonzero vector $v\in T_pM $ there exists 
a vector $u\in T_pM$ such that $\omega_p(v,u)$ is nonzero.
In this case $(M,\omega)$ is called a {\em symplectic manifold}. The symplectic form been non-degenerate,
implies that  the dimension of $M$ must be even.  Unless otherwise stated from now on
we assume that the dimension of $(M,\omega)$ is $2n$.

The first and fundamental  example of a symplectic manifold is $(\mathbb{R}^{2n},\omega_0)$. Here
we take $(x_1,y_1,\ldots, x_n,y_n)$ as coordinates in $\mathbb{R}^{2n}$ and the symplectic form
is defined as
$$
\omega_0:=dx_1\wedge dy_1+\cdots+ dx_n\wedge dy_n.
$$

In the 2-dimensional case, a symplectic form is the same as a volume form.
Hence an oriented surface together with a volume form is an example of a symplectic manifold. 

\begin{example}
Let $S^2$ be the unit sphere in $\mathbb{R}^3$  centered at the origin. Then for
$p\in S^2$ and $u,v\in T_p S^2$ define
$$
\omega_{p}(u,v):= \langle p, u\times v\rangle
$$
where $\langle \cdot, \cdot\rangle$ is the inner product 
and $\times$ is the cross product in  $\mathbb{R}^3$. So defined 
$\omega$ is a non-degenerate $2$-form on the sphere. By dimension reasons, $d\omega=0$;
therefore $\omega$ is a symplectic form on the unit 2-sphere.
\end{example}

\begin{example}
\label{exa:cotang}
Another important class of examples of symplectic manifolds are cotangent bundles $T^*N$
of any smooth manifold $N$. Let $\pi : T^*N\to N$ be the projection map
and $\pi_{*,(q,v^*)}: T_{(q,v^*)} T^*N\to T_qN$ its differential at
$(q,v^*)$. Define the 1-form $\lambda_{\textup{can}}$ on $ T^*N$ at 
$(q,v^*)$ as
$$
\lambda_{\textup{can},(q,v^*)}:= v^* \circ \pi_{*,(q,v^*)}.
$$
Then the canonical symplectic form on $T^*N$ is defined as 
$\omega_\textup{can}:=-d\lambda_{\textup{can}}$. This example is particularly 
important in Classical Mechanics. In fact the
roots of symplectic geometry  go back to Classical Mechanics. For example
see \cite{arnold-math}.

In particular if  $N=\mathbb{R}^n$ with  coordinates $(x_1,\ldots, x_n)$ 
and the fibre $T_x^*\mathbb{R}^n$
with coordinates $(y_1,\ldots, y_n)$, then 
$T^*\mathbb{R}^n$ has coordinates $(x_1,\ldots, x_n,y_1,\ldots, y_n)$.
In this case the 1-form 
defined above takes the form
$$
\lambda_{\textup{can}}=\sum_{j=1}^n y_j dx_j.
$$
Moreover  $\omega_\textup{can}=-d\lambda_{\textup{can}}=\omega_0$;  hence  on 
$T^*\mathbb{R}^n\simeq\mathbb{R}^{2n}$ we get the  
symplectic form defined at the beginning of  this section.
\end{example}

Another source of examples of symplectic manifolds are K\"ahler manifolds. In particular the complex projective space
$(\mathbb{C}P^n,\omega_{FS})$ admits a symplectic form called the {\em Fubini-Study symplectic
form}, which is induced from the Fubini-Study  hermitian metric. That is, if $U_j=\{[z_0:\cdots:z_n] | \,z_j\neq 0\}\subset 
\mathbb{C}P^n$ is a canonical open set, then on $U_j$ the Fubini-Study symplectic 
form is defined as
$$
\omega_{FS}=\frac{i}{2} \partial \overline{\partial} \left( \log  \frac{z_0\overline{z}_0+\cdots +z_n\overline{z}_n}
{z_j\overline{z}_j}
\right).
$$
Furthermore, the symplectic area of the complex line $\mathbb{C}P^1\subset(\mathbb{C}P^n,\omega_{FS})$
is $\pi$.

It is also possible to create new symplectic manifolds from old ones. Cartesian product of symplectic manifolds
is such an example; since this example will be of importance later on 
we explain it in Example \ref{exa:MxM}. Nonetheless there
are more ways to create new symplectic manifolds, such as: symplectic reduction, fibrations and blow ups just
to name a few.

\begin{example}
\label{exa:MxM}
Let $(M,\omega)$ and $(N,\eta)$  a symplectic manifolds
and consider $M\times N$ with projection maps $\pi_M$ and $\pi_N$. Then
$(M\times N, \pi_M^*(\omega)+\pi_N^*(\eta))$ is a symplectic manifold.
In particular for the same symplectic manifold and projections
maps $\pi_j:M\times M\to M$ for $j=1,2$,  we get the symplectic manifold 
$(M\times M, \pi_1^*(\omega)+\pi_2^*(\omega))$.
However it will be more important to consider a different symplectic form
on $M\times M$,
namely $\pi_1^*(\omega)-\pi_2^*(\omega)$.
Below we will see why the minus sign is important  in the second term. 
%If the minus sign is
%replaced by a positive sign, the resulting differential form is still a symplectic form on $M\times M$.
\end{example}
As mentioned above, the standard euclidean symplectic
space $(\mathbb{R}^{2n},\omega_0)$ is the fundamental example of a symplectic manifold.
The reason is that 
locally any symplectic manifold looks like  $(\mathbb{R}^{2n},\omega_0)$.

\begin{theorem}[Darboux]
Let $(M,\omega)$ be a symplectic manifold and $p\in M$. Then there exists
a coordinate chart $(U,\psi)$ about $p$ such
that 
$$
\omega =\psi^*(\omega_0)
$$
on $U$.
\end{theorem}

An important consequence of the above result is that symplectic manifolds
do not have local invariants. 
Thus the techniques and methods used in symplectic geometry
are different from those in Riemannian
geometry.

As mentioned at the Introduction, we aim to give a 
broad overview of Lagrangian Floer homology; which is defined for compact and exact
symplectic manifolds. A symplectic manifold $(M,\omega)$ is called {\em exact} if there exists
a 1-form $\lambda$ such that $\omega=d\lambda$. Each case has its own hypothesis and 
restrictions. In this note we will  cover only  the compact case.
Thus from now on the symplectic manifold $(M,\omega)$ will be assumed to be closed, that
is compact with no boundary. However to illustrate some concepts, some examples will
take place on arbitrary symplectic manifolds.

A {\em Lagrangian submanifold} $L$ of  a symplectic manifold $(M,\omega)$
is an embedded submanifold  $j:L\to M$ of dimension $n$
such that $j^*(\omega)$ is identically zero.
%$L$ 
For example,  the unit circle $S^1$ centered at the origin in 
$(\mathbb{R}^{2},\omega_0)$ is a Lagrangian submanifold. 
More generally, any embedding of $S^1$ into 
 a two-dimensional symplectic manifold is Lagrangian, for dimensional reasons.

On the complex projective space $(\mathbb{C}P^n,\omega_{FS})$, the real projective space
submanifold 
$$
\mathbb{R}P^n=\{[z_0:\cdots:z_n] \in \mathbb{C}P^n | \, z_0,\ldots ,z_n\in \mathbb{R} \}
$$
 is Lagrangian. Another important Lagrangian submanifold of $(\mathbb{C}P^n,\omega_{FS})$ 
 is the {\em Clifford torus} defined as
 $$
\{[z_0:\cdots:z_n] \in \mathbb{C}P^n | \,|z_0|=\cdots =|z_n| \}.
$$
Notice that $
\mathbb{R}P^n$ and the Clifford torus meet in $2^n$ points, namely $[\pm 1:\cdots:\pm 1]$.

In the case of the symplectic manifold $(T^*N,\omega_\textup{can}=-d\lambda_{\textup{can}})$ 
of Example \ref{exa:cotang}, the zero section
and a fiber are examples of  Lagrangian submanifolds. In particular the subspaces
$\{(x_1,0,x_2,0,\ldots 0,x_n,0)\}$ and $\{(0,y_1,0,y_2,0,\ldots, y_{n-1},0,y_n)\}$ of
$(\mathbb{R}^{2n},\omega_0)$ are Lagrangian submanifolds. There is also another significant class
of Lagrangian submanifolds of $(T^*N,-d\lambda_\textup{can})$.  Let $\sigma$ be a 1-form on $N$
and consider it as a section  $\sigma: N\to T^*N$; that is $\sigma(p)=\sigma_p(\cdot)$. Hence
$\sigma$ embeds $N$ into $T^*N$ and  following the definition of the canonical 1-form we get that 
$$
\sigma^*(\lambda_{\textup{can}})=\sigma.
$$
Hence the graph of a 1-form is a Lagrangian submanifold of $(T^*N,-d\lambda_\textup{can})$ 
if and only if the 1-form is closed.

In the case of the symplectic manifold 
$(M\times M, \pi_1^*(\omega)-\pi_2^*(\omega))$ of Example \ref{exa:MxM}, the diagonal
$\Delta=\{(x,x)|\, x\in M\}$ is a Lagrangian submanifold. Notice that the minus sign in the
second term of the symplectic form in fundamental to guarantee that $\Delta$ is Lagrangian.
Below in Example \ref{exa:graph} 
we will exhibit 
more Lagrangian submanifolds of $(M\times M, \pi_1^*(\omega)-\pi_2^*(\omega))$ that
are related to symplectic diffeomorphisms of $(M,\omega)$; the diagonal $\Delta$ is 
the graph of the identity diffeomorphism.

The relevance  of the symplectic manifold $(T^*N, -d\lambda_{\textup{can}})$ of Example 
\ref{exa:cotang}, is that it is a symplectic model of a neighborhood whenever
$N\subset M$ is a Lagrangian submanigfold, regardless of the symplectic manifod
$(M,\omega)$.
That is, if $N$ is a Lagrangian submanifold of  $(M,\omega)$ 
then a tubular neighborhood of it can be identified, in a suplectic way, with
a neighborhood of the zero section of $(T^*N, -d\lambda_{\textup{can}})$. 
That is, there is an analog of  Darboux's Theorem for Lagrangian submanifolds, where the standard 
symplectic euclidean space
$(\mathbb{R}^{2n},\omega_0)$
is replaced by the standard symplectic cotangent bundle 
$(T^*N, -d\lambda_{\textup{can}})$. Recall from above that in the  case $N=\mathbb{R}^n$, 
we showed that $(T^*\mathbb{R}^n, -d\lambda_{\textup{can}})$ agrees with $(\mathbb{R}^{2n},\omega_0)$.

\begin{theorem}[Weinstein]
Let $(M,\omega)$ be a symplectic manifold and $L$ a Lagrangian submanifold. Then there exists
a neighborhood $U$ of $L$ and a neighborhood $V$ of $L_0$, the zero section of  
$(T^*L, -d\lambda_{\textup{can}})$, that are diffeomorphic by $\psi:U\to V$ such that 
$$
\omega =\psi^*(-d\lambda_{\textup{can}})
$$
and $\psi(L)=L_0.$
\end{theorem}

Using  the fact that locally any symplectic manifold is equal to $(\mathbb{R}^{2n},\omega_0)$,
is possible to  show that there are many Lagrangian submanifolds in any  given symplectic manifold 
$(M,\omega)$.
The circle
 $S^1$ is Lagrangian submanifold of $(\mathbb{R}^{2},\omega_0)$. Moreover we can make the
radius arbitrary small, say $\epsilon>0$,  and still $S^1(\epsilon)$ is a Lagrangian 
submanifold. Taking $n$ copies
 of this example, it follows that the $n$-dimensional 
  $\epsilon$-torus  $S^1(\epsilon)\times \cdots\times S^1(\epsilon)$ is a Lagrangian submanifold
of $(\mathbb{R}^{2},\omega_0)\times \cdots\times(\mathbb{R}^{2},\omega_0) =
(\mathbb{R}^{2n},\omega_0)$. 
Thus for a given  symplectic manifold $(M,\omega)$ and
$\epsilon $  small enough,  by Darboux's Theorem we have that 
the $n$-dimensional   $\epsilon$-torus is a Lagrangian submanifold of $(M,\omega)$.
% of , an arbitrary symplectic manifold.

In the particular case of $(\mathbb{R}^{2},\omega_0)$, we have that the two-dimensional 
torus $S^1\times S^1$ is a Lagrangian submanifold. Furthermore, the torus is the only oriented
surface that can be embedded as a Lagrangian submanifold in  $(\mathbb{R}^{2},\omega_0)$.

%%%%%%%%%%%%%%%%%%%%%%%%%%%%%%%%%%%%%%%%%%%%%%
%%%%%%%%%%%%%%%%%%%%%%%%%%%%%%%%%%%%%%%%%%%%%%
\section{Symplectic and Hamiltonian Diffeomorphisms}
%%%%%%%%%%%%%%%%%%%%%%%%%%%%%%%%%%%%%%%%%%%%%%
%%%%%%%%%%%%%%%%%%%%%%%%%%%%%%%%%%%%%%%%%%%%%%

There are two types of symmetries associated to a symplectic manifold. Recall that
we assumed that that symplectic manifold is closed.
%In order to avoid
%some minor issues, we will assume that the manifold $M$ is closed; that is compact with
%no boundary. 
In the non-compact case, one has to consider diffeomorphisms with compact support. 
A diffeomorphism $\phi:(M,\omega)\to (M,\omega)$ is said to be a  {\em symplectic diffeomorphism}
if $\phi^*(\omega)=\omega$. The set of symplectic diffeomorphisms of $(M,\omega)$
forms a group under composition 
and is denoted by $\textup{Symp}(M,\omega)$. 
In fact the group of symplectic diffeomorphisms is an infinite dimensional space,
its Lie algebra consists of vector fields $X$ such that the
1-form $\omega(X, \cdot)$ is closed.

Among the group of symplectic diffeomorphisms
we have the second type of symmetries, called Hamiltonian diffeomorphisms.
A symplectic diffeomorphism $\phi$ is called  {\em  Hamiltonian diffeomorphism} if there exists
a path of symplectic diffeomorphisms
 $\{\phi_t\}_{0\leq t\leq 1}$ and a
smooth function $H:[0,1]\times M\to \mathbb{R}$,
such
that $\phi_0=1_M$, $\phi_1=\phi$, and if  $X_t$ is the time-dependent vector field
induced by the equation
$$
\frac{d}{dt} \phi_t=X_t\circ \phi_t,
$$
then $\omega(X_t,\cdot)=dH_t$. 
The set of Hamiltonian diffeomorphisms is a group under composition and is denoted by 
$\textup{Ham}(M,\omega)$.
As in the symplectic case, $\textup{Ham}(M,\omega)$ is an infinite dimensional space and its
Lie algebra consists of vector fields $X$ such that the 1-form $\omega(X,\cdot)$ is exact.

A Hamiltonian diffeomorphisms $\phi$ is  called {\em autonomous}, if there  exists a  path $\{\phi_t\}$,
 as in the definition of Hamiltonian diffeomorphism,  such that 
$X_t$  is independent of $t$. In other words autonomous Hamiltonian
diffeomorphisms are the image of the exponential map of Hamiltonian vector fields. 
%$X$ on $(M,\omega)$ such that $\omega(X,\cdot)$ is an exact 1-form.
Alternatively,  the group $\textup{Ham}(M,\omega)$ can be described as the group generated
by autonomous Hamiltonian diffeomorphisms \cite{banyaga-surla}.

Not only  is $\textup{Ham}(M,\omega)$ a subset 
$\textup{Symp}(M,\omega)$, as the definition suggest;  the group of Hamiltonian diffeomorphisms
is a normal subgroup of the group of symplectic diffeomorphisms. 
As we explain below in most cases is a proper subgroup.
Among other properties of  $\textup{Ham}(M,\omega)$ is that it is connected with respect to the
$C^\infty$-topology; $\textup{Symp}(M,\omega)$ does not have to be connected. Further if
$\textup{Symp}_0(M,\omega)$ is the connected component of the group of symplectic diffeomorphisms
that contains the identity map and  $H^1(M,\mathbb{R})=0$, then 
\begin{eqnarray*}
%\label{e:subgroup}
\textup{Ham}(M,\omega)=\textup{Symp}_0(M,\omega).
\end{eqnarray*}
%In fact, if $H^1(M,\mathbb{R})=0$ then the above inclusion is an equality. Thus f
For example
$\textup{Ham}(\mathbb{C}P^n,\omega_{\textup{FS}})=\textup{Symp}_0(\mathbb{C}P^n,
\omega_{\textup{FS}})$ for $n\geq 1.$ 
However if $H^1(M,\mathbb{R})\neq 0$, 
$\textup{Ham}(M,\omega)$ is properly contained in
$\textup{Symp}_0(M,\omega).$
Below we will see an example where
$\textup{Ham}(M,\omega)$
 is a proper subgroup of $\textup{Symp}_0(M,\omega)$.

An important remark about a Hamiltonian diffeomorphism is that  its fixed point set 
of  is non-empty. 
Recall that we are assuming that $(M,\omega)$ is closed; in the non-compact case the assertion
is false. For instance, on $(\mathbb{R}^{2n},\omega_0)$ a translation map
is a Hamiltonian diffeomorphism that is fixed-point free.
As for the case of compact symplectic manifolds 
it is straightforward to justify that the fixed point set is non-empty in the  case of 
autonomous Hamiltonians.  For if $\phi$ is an  autonomous Hamiltonian, 
then there exists $\{\phi_t\}$ such that
\begin{eqnarray}
\label{e:relation}
\frac{d}{dt} \phi_t=X\circ \phi_t,  \qquad \textup{and}  \qquad \omega(X,\cdot)=dH.
\end{eqnarray}
Since the manifold is assumed to be compact, then the set of critical points
of $H:M\to \mathbb{R}$ is non-empty. But $\omega$ is non-degenerate, hence by
 Eqs. (\ref{e:relation}) the set of critical points of $H$
coincides with the zero set of $X$. If $X$ vanishes at $p$ it follows by Eqs. (\ref{e:relation}) that  
$p$ is a fixed point of the flow
$\{\phi_t\}$; in particular $p$ is a fixed point of $\phi_1=\phi.$

The fact that Hamiltonian diffeomorphisms on compact manifold always have fixed
points is not shared by symplectic diffeomorphisms that are not Hamiltonian.

\begin{example}
\label{exa:torusmap}
Consider the flat two-dimensional torus $(\mathbb{T}^2=\mathbb{R}^2/\mathbb{Z}^2,
\omega=dx\wedge dy)$.
For a fix $\alpha\in (0,1)$, the translation  map
$$
\phi^\alpha[x,y]:=[x+\alpha,y]
$$
preserves the area and hence is symplectic diffeomorphism. Notice that since $\alpha\neq 0,1$, the map
$\phi^\alpha$ has no fixed points. Hence $
\phi^\alpha$ is not a Hamiltonian diffeomorphism for any $\alpha\in (0,1)$. Moreover $\phi^\alpha$
lies in the identity component of the group of symplectic diffeomorphism.
Therefore
$\textup{Ham}(\mathbb{T}^2,\omega)$ is a proper subgroup of $\textup{Symp}_0(\mathbb{T}^2,\omega)$.
\end{example}

As mentioned above, symplectic diffeomorphisms give rise to Lagrangian submanifolds. In the next
example we show how this is done and highlight the importance of this example 
in the study of fixed points of Hamiltonian diffeomorphisms.

\begin{example}
\label{exa:graph}
Let $\phi:(M,\omega)\to (M,\omega)$ be a symplectic diffeomorphism,
thus $\phi^*(\omega)=\omega$.
Then the graph of $\phi$ is an embedded submanifold of dimension $2n$ in $M\times M$,
that is $j:M\to M\times M$ is given by $j(x)=(x,\phi(x))$ and its image 
is the graph of $\phi$,
$$
\textup{graph}(\phi):=\{ (x, \phi(x)) | \,x\in M  \}.
$$
Furthermore the graph of $\phi$ is a Lagrangian submanifold of
$(M\times M, \pi_1^*(\omega)-\pi_2^*(\omega))$; for
\begin{eqnarray*}
j^*(    \pi_1^*(\omega)-\pi_2^*(\omega)) =   \omega -\phi^*(\omega)=0.
\end{eqnarray*}
The above computation shows the relevance of the minus sign that appears in the symplectic form of 
$(M\times M, \pi_1^*(\omega)-\pi_2^*(\omega))$ in Example \ref{exa:MxM}.
In the case when $\phi:M\to M$ is a Hamiltonian diffeomorphism, when know that
the fixed point set is non empty; further this set is in one-to-one correspondence
with the intersection points of  $\textup{graph}(\phi)$ with the diagonal $\Delta$. As pointed out
above, there are symplectic diffeomorphisms $\phi$, such that $\textup{graph}(\phi)$ and $\Delta$
have no points in common. For instance $\phi^\alpha$ of Example
\ref{exa:torusmap}.
\end{example}

Notice that for any Hamiltonian diffeomorphism $\phi:(M,\omega)\to (M,\omega)$ and 
Lagrangian submanifold $L\subset (M,\omega)$,
$\phi(L)$ is again a Lagrangian submanifold.
An important fact that will be useful in the context of Lagrangian Floer homology is the following.
A Lagrangian submanifold $L\subset (M,\omega)$ is called {\em non-displaceable} if for 
every Hamiltonian diffeomorphisms
$\phi:(M,\omega)\to (M,\omega)$, the Lagrangian submanifolds $L$ and $\phi(L)$ have points 
in common. Otherwise, $L$
is   called {\em displaceable}. 
Hence we are considering the intersection of two particular Lagrangian submanifold,
$L$ and $\phi(L)$. This is part of the phenomenon that Lagrangian Floer homology attempts 
to answer, intersection or non-intersection of Lagrangian submanifolds.

Consider the two-dimensional sphere $(S^2,\omega)$ with any area form, 
let us try to understand the intersection of a particular pair of Lagrangian submanifolds.
Consider the Lagrangian submanifold $L$ to be any circle
that lies entirely in a hemisphere. 
Thus there always exists a rotation
$\phi:(S^2,\omega)\to (S^2,\omega)$, which is in fact a Hamiltonian diffeomorphism of the 2-sphere, such 
that $L$ and $\phi(L)$ have no points in common. That is $L$ is displaceable.  
Now consider the  case when $L$
is such that  both components $U$ and $ V$ of $S^2\setminus L$ have equal area. Recall that 
any Hamiltonian diffeomorphism preserves area;  hence for any  Hamiltonian diffeomorphism $\phi$,
$\phi(U)=U$  or $\phi(U)\cap V$ is not empty. 
Then for any Hamiltonian diffeomorphism $\phi$, we have that  $L\cap\phi(L)$ is non-empty
if the Lagrangian submanifold is such that the two components of $S^2\setminus L$ have equal area.
\begin{comment}
In general, a Lagrangian submanifold of  $(M,\omega)$ 
 is called {\em non-displaceable}
if for every $\phi$ Hamiltonian diffeomorphisms the set $L\cap \phi(L)$ is non-empty.
Notice that
$ \phi(L)$ is  a  
Lagrangian submanifold of  $(M,\omega)$  for any Hamiltonian diffeomorphism $\phi$.
Thus in the case of $(S^2,\omega)$, 
\end{comment}
Hence the non-displaceable Lagrangian are precisely the
embedded circles that split the sphere in two pieces of equal area.
In fact one of the current problems in symplectic geometry is to determine which Lagrangian
submanifolds are non-displaceable or displaceable.

The higher dimensional analog of the above example, for the case
when $L$ splits the sphere in two parts of equal area, is the Lagrangian submanifold 
$\mathbb{R}P^n$ in  $(\mathbb{C}P^n, \omega_{FS})$. One of the triumphs of 
Lagrangian Floer homology is the proof that  
$\mathbb{R}P^n$ is non-displaceable. This result was proved by  Y.-G. Oh in \cite{oh-floercoho}; where
he defined Lagrangian Floer homology  for monotone Lagrangian submanifolds.

Now we go back to the case of Lagrangian submanifolds induced by symplectic
diffeomorphisms as in  Example \ref{exa:graph}. Hence let $\phi:(M,\omega)\to (M,\omega)$ 
be a symplectic
diffeomorphism
and $\textup{graph}(\phi)  \subset (M\times M, \pi_1^*(\omega)-\pi_2^*(\omega))$
which is a Lagrangian submanifold.
Note that in this example
 the Lagrangian submanifold $\textup{graph}(\phi)$  it is actually the image of the 
Lagrangian $\Delta$
under the symplectic diffeomorphisms $1\times\phi$ 
of $(M\times M, \pi_1^*(\omega)-\pi_2^*(\omega))$. That is
$$
 \textup{graph}(\phi)= (1\times \phi)(\Delta).
$$
In fact when $\phi:(M,\omega)\to (M,\omega)$ is a Hamiltonian diffeomorphism we know that $ (1\times \phi)(\Delta)\cap \Delta$
is non empty. The intersection points are in one-to-one correspondence with the fixed
points of $\phi$. In this case $ (1\times \phi)$ is also a Hamiltonian diffeomorphism
of $(M\times M, \pi_1^*(\omega)-\pi_2^*(\omega))$.
Lagrangian Floer homology gives a stronger result, it shows that $\Delta$
is non displaceable, that is $\Phi(\Delta)\cap \Delta\neq \emptyset$ for {\bf any}
Hamiltonian diffeomorphisms ${\Phi}$
of  $(M\times M, \pi_1^*(\omega)-\pi_2^*(\omega))$, 
not necessarily those induced from Hamiltonians of $(M,\omega).$
 Moreover it gives a lower
bound on the cardinality of $\Phi(\Delta)\cap \Delta$ under some non degeneracy conditions
of $\Phi$.
That is it solves the Arnol'd Conjecture. % that we will  discuss below.

\begin{comment}
That is, in the case of Example \ref{exa:torusmap} we have that for any $\alpha\in (0,1)$, the intersection 
$\Delta\cap (1\times\phi^\alpha)(\Delta)$ is empty.
However 
$\Delta$ is  non-displaceable Lagrangian submanifold of 
$(\mathbb{T}^2\times \mathbb{T}^2,
\pi_1^*\omega-\pi_2^*\omega )$ since any Hamiltonian diffeomorphisms. However if we consider $(M,\omega)$ a 
closed symplectic manifold,  then
$$
\Delta\cap (1\times\phi)(\Delta)\neq \emptyset
$$
and every Hamiltonian diffeomorphism $\phi$. That is to say, 
$\Delta$ is a non-displaceable Lagrangian submanifold of $(M\times M, \pi_1^*(\omega)-\pi_2^*(\omega))$.
\end{comment}

The problem of estimating the number of fixed points of a Hamiltonian
diffeomorphism, is a particular case of the wider problem of estimating the number 
of intersection points of two Lagrangian submanifolds.
In broad terms, that  is the objective of Lagrangian Floer
homology.

As seen in the definition, Hamiltonian diffeomorphisms have a strong connection with 
 smooth functions. A manifestation  of this connection was the nice link between the fact that a 
smooth function on a closed manifold admits critical points;  and the fact the on a closed symplectic
manifold  the fixed point set of a Hamiltonian diffeomorphism is non-empty.

In 1965, V.  Arnol'd \cite{arnold-sur} conjectured an analog result of
Theorem \ref{t:morseine}, but for the case
of Hamiltonian diffeomorphisms on closed symplectic manifolds 
 instead  of Morse functions on arbitrary manifolds.
See also \cite[Appendix 9]{arnold-math}. His motivation was
the Poincar\'e-Birkhoff annulus theorem: An area preserving diffeomorphism of the annulus 
such that the boundary circles are turned in different  directions must have at least two fixed points.
A fixed point $p\in M$ of a Hamiltonian diffeomorphism $\phi$ is said to be {\em non-degenerate} if
$1$ is not an eigenvalue of the  the linear map $\phi_{*,p}: T_pM\to T_pM$. Note that
non-degenerate fixed points are isolated, and in the case of a closed symplectic manifold there are
a finite number of them.

\begin{conjecture}[Arnol'd]
Let $(M,\omega) $ be a closed symplectic manifold and $\phi$ a Hamiltonian
diffeomorphism such that all of  its fixed points are non-degenerate. Then 
$$
\#\{p\in M |  \phi(p)=p\} \geq   \sum_{j=0}^{2n} \textup{Rank }  H_j(M,\mathbb{R}).
$$
\end{conjecture}

For two-dimensional symplectic manifolds,  the conjecture was proved by Y. Eliashberg  \cite{eliashberg-atheorem}; 
in \cite{conley-zehnder-thebir}
C. C. Conley and E. Zehnder proved the conjecture for the symplectic torus manifold with the standard symplectic form;
and for the complex projective space with the Fubini-Study symplectic form  the conjecture was
proved by
B. Fortune and A. Weinstein in \cite{fortune-weinstein-asymp}.
The real break through in solving Arnold's conjecture was made by A. Floer in
\cite{floer-morse}. 

 In a series of papers  A. Floer developed a homological theory based
on holomorphic techniques, which were introduced by   M. Gromov \cite{gromov-psudo}, and the new
approach to Morse theory developed by E. Witten \cite{witten-super}.
Under some assumption on the symplectic manifold A. Floer developed Hamiltonian
Floer homology, using holomorphic cylinders, in order to find a lower bound to the number of fixed points
of a Hamiltonian diffeomorphism.  Then he generalized 
this approach to develop Lagrangian Floer homology, now using holomorphic stripes,
in order to determine the minimum number of intersection
points of a particular pair of Lagrangian submanifolds.

The Arnold's conjecture has been proved for arbitrary symplectic manifolds.
Some reference for the proof of the conjecture, sometimes under some restrictions and
others in full generality are:
 K. Fukaya and K. Ono \cite{fukaya-ono-arnold}, 
 H. Hofer and D. Salamon \cite{hofer-salamon-floerhomo}, 
G.  Liu and G. Tian  \cite{liu-tian-floerhomo}, K. Ono \cite{ono-onthe}, 
Y.-G. Oh \cite{oh-floercoho},
 Y. Ruan \cite{ruan-virtual}.
    
Some of the techniques introduced in \cite{fukaya-ono-arnold} on the proof 
of the Arnold's conjecture have been re-evaluated.
For instance the
Kuranishi structure on the moduli space of holomorphic strips 
$ \mathcal{M}_J(p,q,L_0,L_1)$, that will be defined  in the next section, as well
as its virtual fundamental class. 
Recently J. Pardon \cite{pardon-analgebraic} has given an alternative approach 
to this problem using techniques from homological algebra.
There is also a series of articles by D. Mcduff and K.  Wehrheim, \cite{mcduff-notes}
\cite{mw-1}, \cite{mw-2} and \cite{mw-3}, where they treat this problem using 
tools from analysis.

%\begin{proposition} %%%%Ver Audin pag.137
%A non-degenerate critical point of $H$ corresponds to a non-degenerate
%critical point of $\phi$. Conversely,  if $H$ is $C^2$-small then the two concepts
%agree.
%\end{proposition}

%%%%%%%%%%%%%%%%%%%%%%%%%%%%%%%%%%%%
%%%%%%%%%%%%%%%%%%%%%%%%%%%%%%%%%%%%
\section{Lagrangian Floer Homology}
%%%%%%%%%%%%%%%%%%%%%%%%%%%%%%%%%%%%
%%%%%%%%%%%%%%%%%%%%%%%%%%%%%%%%%%%%

The construction of Lagrangian Floer homology emulates to a certain extent 
the construction of  Morse homology described above.
The manifold in consideration to define it is a certain 
space of trajectories which is infinite dimensional; and the function defined on it
is a certain action functional. 
The critical points turn out to be constant
trajectories, that give rise to the differential complex used to define
Lagrangian Floer homology. 
Is important to point out while the construction of Lagrangian Floer homology 
follows the spirit of the construction of Morse homology, new complications 
emerge that were not present before. 
Just to have an idea of this, it suffices to say that Lagrangian Floer homology 
is not always defined do to the fact that the square of the differential map is not always
equal to  zero.

Let $L_0$ and $L_1$ be two compact Lagrangian submanifolds in $(M,\omega)$ 
that intersect transversally. 
Consider the space of smooth trajectories from $L_0$ to $L_1$,
$$
\mathcal{P}(L_0,L_1):=\{ \gamma:[0,1]\to M  |   \gamma \textup{ is smooth, } 
\gamma(0)\in L_0 \textup{ and }\gamma(1)\in L_1 \}.
$$
endowed with the $C^\infty$-topology.
Notice that the constant paths  in $
\mathcal{P}(L_0,L_1)$ are the ones that correspond to the intersection points $L_0\cap L_1$.
From now on we write $\mathcal{P}$ for $
\mathcal{P}(L_0,L_1)$.

The space $\mathcal{P}$ is not necessarily connected, thus we 
fix $\hat\gamma$ in $\mathcal{P}$ and consider  the component that
contains $\hat\gamma$,  which we denoted by  $\mathcal{P}(\hat\gamma)$. 
Relative to  $\hat\gamma$  consider the 
universal covering space
$\widetilde{\mathcal{P}}(\hat\gamma)$ of ${\mathcal{P}}(\hat\gamma)$. Elements of
$\widetilde{\mathcal{P}}(\hat\gamma)$ are denoted by $[\gamma, w]$ where
$w$ is a smooth path in $\mathcal{P}(\hat\gamma)$ from   $\hat\gamma$ to   $\gamma$. That is
$w:[0,1]\times [0,1] \to M$ is a smooth  map such that 
$w(s,\cdot)\in \mathcal{P}(\hat\gamma)$ for all $s\in [0,1]$,
$w(0,\cdot)= \hat\gamma$ and
 $w(1,\cdot)=\gamma$.

The space $\widetilde{\mathcal{P}}(\hat\gamma)$ is not the right space to define the
action functional.  The right space is the Novikov covering of ${\mathcal{P}}(\hat\gamma)$,
which is defined by an equivalence relation on $\widetilde{\mathcal{P}}(\hat\gamma)$.
For the sake of making the exposition less technical we are not going to define
the Novikov covering of ${\mathcal{P}}(\hat\gamma)$, instead we are going to impose
strong assumptions on the symplectic manifold $(M,\omega)$ and the pair of 
Lagrangian submanifolds $L_0$ and $L_1$
in order to define the action functional on $\widetilde{\mathcal{P}}(\hat\gamma)$
and carry out a similar procedure as in Morse theory.
Thus from now on we assume that the symplectic manifold $(M,\omega)$ is such
that
$$
\int_{S^2} f^*\omega=0
$$
for every $[f]\in \pi_2(M)$, where $f$ is a smooth representative.  
A symplectic manifold that satisfies this
condition is said to be {\em symplectically aspherical}. 
The symplectic tori
$(\mathbb{T}^{2n},\omega)$, for   $n\geq 1$,  are examples symplectically aspherical
since $\pi_2(\mathbb{T}^{2n})$ is trivial. However there are plenty of symplectically aspherical
manifolds with non trivial $\pi_2$, even in dimension four \cite{gompf-symplectically}.
As we will see in the next paragraph, if $(M,\omega)$ is symplectically aspherical then the action
functional is  well defined on the covering
$\widetilde{\mathcal{P}}  (L_0,L_1;\hat\gamma)$. This hypothesis on $(M,\omega)$ is also useful
since it rules out the appearance of bubbles; that is holomorphic spheres attached to 
holomorphic strips. See for example \cite{}. As mentioned above, Lagrangian Floer homology
 is defined on more generally symplectic manifolds, even in the presence of bubbles.

Further  we also assume that $L_0$ and $L_1$ are simply connected.
Then the  action functional $\mathcal{A}:  \widetilde{\mathcal{P}}  (L_0,L_1;\hat\gamma)\to \mathbb{R }$ 
is defined as
$$
\mathcal{A}([\gamma, w])=\int_{[0,1]\times [0,1] } w^*\omega,
$$
that is, the symplectic area of $w([0,1]\times [0,1]  )\subset (M,\omega)$. 
To see that $\mathcal{A}$ is well defined, let
$(\gamma,w)$ and  $(\gamma,w^\prime)$ represent the same class.
%For  two representatives of the same
%class $[\gamma,w]=[\gamma,w^\prime]$, 
Then we have a map defined on the cylinder,
$\overline{w}\#w^\prime:S^1\times [0,1]\to M$
where $\overline{w}\#w^\prime(s,0)$ is a loop in $L_0$ and 
$\overline{w}\#w^\prime(s,1)$ is a loop in $L_1$. Since the Lagrangian submanifolds are
assumed to be simply connected,
there exists a 2-disk contained in $L_0$ whose boundary is the loop $\overline{w}\#w^\prime(s,0)$. 
This observation also applies to $L_1$. That is, we added the caps to the cylinder  
to obtain topological  a  2-sphere. Since 
$(M,\omega)$
is symplectically aspherical, the symplectic area of the 2-sphere is zero. Notice also that the symplectic
area of each cap is zero, since the symplectic form is identically zero on Lagrangian submanifolds.
Thus the symplectic area of the cylinder
$\overline{w}\#w^\prime:S^1\times [0,1]\to M$ is equal to zero. Then we have that 
\begin{eqnarray*}
0&=&\int_{S^1\times [0,1] }\overline{w}\#w^*\omega\,\,=
\int_{[0,1]\times [0,1] } \overline{w}^*\omega\,\,+\int_{[0,1]\times [0,1] } w^*\omega\\
&=&-\int_{[0,1]\times [0,1] } {w}^*\omega\,\,+\int_{[0,1]\times [0,1] } w^*\omega
\end{eqnarray*}
Hence it follows  that the action functional is well defined on
the covering $\widetilde{\mathcal{P}}  (L_0,L_1;\hat\gamma)$.
In local coordinates the action functional
takes the form
$$
\mathcal{A}([\gamma, w])=\int_0^1 \int_0^1 \omega\left( \frac{\partial w}{\partial s} ,
\frac{\partial w}{\partial t}\right) ds \, dt.
$$

Lagrangian Floer homology is defined  emulating the way Morse homology is defined.
The finite dimensional manifold in Morse homology is replaced by the infinite dimensional
space $\widetilde{\mathcal{P}}  (L_0,L_1;\hat\gamma)$. And the function into consideration
is the action functional. However the analytical difficulties in this setting are more complex than in the
Morse scenario.

To follow the path of Morse homology we need to define a Riemannian structure on 
$\widetilde{\mathcal{P}}  (L_0,L_1;\hat\gamma)$.  Denote by $\mathcal{J}$ the space of almost
complex structures on $(M,\omega)$.
Recall that  an almost complex structure $J\in \mathcal{J}$
on $(M,\omega)$ is said to be $\omega${\em-compatible}  if for every nonzero
vector $v$, 
$$
\omega (v ,J v) >0 \quad  \textup{ and }   \quad  \omega (J\cdot ,J \cdot) = \omega (\cdot ,\cdot).  
$$
For a $\omega$-compatible almost complex structure $J$, we have that
$$
g_J(\cdot ,\cdot):=\omega (\cdot ,J \cdot)
$$
defines a Riemannian metric on $(M,\omega)$. 
%Let $\mathcal{J}_\omega$ be the set of compatible
%almost complex structures.
%is a Riemannian metric on $M$.  
Is important to note that 
compatible almost complex structures exist in abundance
on any symplectic manifold. Let $J=\{J_t\}_{0\leq t\leq 1}$ be a smooth family of 
$\omega$-compatible almost complex
structures on $(M,\omega)$; hence we have  $\{g_t\}_{0\leq t\leq 1}$ a smooth family 
 of Riemannian metrics. Then on 
$\widetilde{\mathcal{P}}(\hat\gamma)$
we define a Riemannian metric associated to $\{g_t\}_{0\leq t\leq 1}$ as
$$ 
\langle\!\!\langle \xi_1,\xi_2\rangle\!\!\rangle: =\int_0^1  g_t(\xi_1(t), \xi_2(t)) dt
$$
for $\xi_1,\xi_2$ in 
$T\widetilde{\mathcal{P}}(\hat\gamma)$.
As in the Morse theory case, we compute the gradient of $\mathcal{A}:  
\widetilde{\mathcal{P}}  (\hat\gamma)\to \mathbb{R }$
with respect to $\langle \!\!\langle\cdot,
\cdot\rangle\!\!\rangle$,
\begin{eqnarray*}
d\mathcal{A}_{([\gamma, w])}(\xi) 
&=& \int_0^1  \omega  \left( \frac{\partial \gamma}{\partial t} ,\xi(t)\right)  dt \\
&=& \int_0^1  \omega  \left( \frac{\partial \gamma}{\partial t} ,J_t(-J_t\xi(t))\right)  dt \\
&=& \int_0^1  g_t  \left( \frac{\partial \gamma}{\partial t} ,-J_t\xi(t)\right)  dt \\
&=&  \left\langle\!\!\!\!\left\langle \frac{\partial \gamma}{\partial t} ,-J_t\xi   \right\rangle\!\!\!\!\right\rangle\\
 &=&\left\langle\!\!\!\!\left\langle J_t\frac{\partial \gamma}{\partial t} ,\xi   \right\rangle\!\!\!\!\right\rangle.
\end{eqnarray*}
That is,
 $$
\textup{grad}\mathcal{A}([\gamma,w])=J_t\frac{\partial \gamma}{\partial t}.
$$

Since $J_t$ is an  automorphism of $TM$ for each $t$, 
the gradient of $\mathcal{A}$ vanishes at $[\gamma,w]$ if and only if  $\gamma:[0,1]\to 
(M,\omega)$ is a constant path.
Thus the critical points of the action functional $\mathcal{A}$ are of the form $[\gamma,w]$  where 
$\gamma$ is a  constant path, corresponding to  an  intersection point
of  $L_0$ with $\ L_1$.

%Let $p$ and $q$ in  $L_0\cap L_1\subset \widetilde{\mathcal{P}}(\hat\gamma)$ be critical points  
%of $\mathcal{A}$, that is we write $p$ for $[p,w]$ where 
%$w:[0,1]\times [0,1] \to M$ is such that 
% $w(1,t)=p$ for all $t$. Similarly for the other intersection point $q$.

 As in the Morse theory case, a  flow line of  $-\textup{grad}\mathcal{A}$ connecting $p $ to $q$
is a  smooth function   $u:\mathbb{R}\to  \widetilde{\mathcal{P}}(\hat\gamma)$
such that 
\begin{eqnarray}
\label{e:flowhol}
\frac{d u}{d s}= -\textup{grad}\mathcal{A},
\quad  \lim_{s\to -\infty} u(s)=p,  \quad \hbox{and}\quad \lim_{s\to +\infty} u(s)=q.
\end{eqnarray}
Unwrapping this,  is better to write $u:\mathbb{R}\times [0,1]\to  (M,\omega)$ and 
the first equation of Eqs. (\ref{e:flowhol}) as
\begin{eqnarray}
\label{e:hol}
\frac{\partial u}{\partial s}+J_t\frac{\partial u}{\partial t}=0.
\end{eqnarray}
Note that Eq. (\ref{e:hol}) is the Cauchy-Riemann
equation, $\overline{\partial}_J(u)=0$, with respecto to the $\omega$-compatible almost complex,
$\{J_t\}$,
structure on $(M,\omega)$. A smooth map $u:\mathbb{R}\times [0,1]\to  (M,\omega)$
is said to be  a
{\em
  $J$-holomorphic
strip} in $(M,\omega)$ if $\overline{\partial}_J(u)=0$. 
 Then the space of connecting flow lines (actually $J$-holomorphic strips in $(M,\omega,J)$)
that connect $p$ with $q$ is defined as
 as
\begin{eqnarray*} 
\hat{ \mathcal{M}}_J(p,q,L_0,L_1):= \{u:\mathbb{R}\times [0,1]\to (M,\omega)&|& u \textup{ is smooth,  satisfies
Eqs. (\ref{e:flowhol})}\\
& & \textup{and }u(s,\cdot)\in \mathcal{P}(L_0,L_1) 
\}.
\end{eqnarray*}
In the case when $(M,\omega)$ is non compact but exact, and the Lagrangian 
submanifolds are still compact,
 one imposes an additional condition on the flow lines. That is to say, in addition 
to Eqs. (\ref{e:flowhol}) the map $u:\mathbb{R}\times [0,1]\to (M,\omega)$ is  
 required to have finite  energy,
\begin{eqnarray}
\label{e:fener}
\int_{\mathbb{R}\times [0,1]}u^*\omega < \infty.
\end{eqnarray}
In the  case of a compact symplectic manifold,
 a holomorphic strip $u$  with $u(s,\cdot)\in \mathcal{P}(L_0,L_1)$
has finite energy,
Eq. (\ref{e:fener}), if and only if satisfies the limit conditions
$$
\lim_{s\to -\infty} u(s)=p,  \quad \hbox{and}\quad \lim_{s\to +\infty} u(s)=q
$$
for some $p,q\in L_0\cap L_1$. For the details see J. Robbin and D. Salamon \cite{robbin-salamon}.

\begin{comment}
To define the space of connecting flow lines (actually $J$-holomorphic strips in $M$), we must add 
one more condition. Namely we require that each solution $u$ of (\ref{e:hol}) must have finite
energy,
\begin{eqnarray}
\label{e:fener}
\int_{\mathbb{R}\times [0,1]}u^*\omega < \infty.
\end{eqnarray}
Then  for $p, q\in L_0\cap L_1$, we denote by $\mathcal{M}(p,q)$ the space of
 flow lines that connect $p$ with $q$ as
 $$
 \mathcal{M}(p,q):= \{u:\mathbb{R}\times [0,1]\to (M,\omega)| u \textup{ is smooth and satisfies (\ref{e:flowhol})  
 and (\ref{e:fener}) }\}.
 $$
\end{comment}

Note that the strip $\mathbb{R}\times [0,1]i\subset \mathbb{C}$ is conformally equivalent 
with the closed unit disk $D^2\subset \mathbb{C}$ minus two points on the boundary. 
Thus sometimes $u$ is also referred as a holomorphic disk.

For $p,q\in L_0 \cap L_1$, let $C^\infty(\mathbb{R}\times [0,1], M; L_0,L_1)$ be the set of smooth maps
$u:\mathbb{R}\times [0,1]\to M$ with the limit behavior as in (\ref{e:flowhol}). Thus we have a bundle map
$$
\cup_u \,C^\infty(u^*(TM))\to C^\infty(\mathbb{R}\times [0,1], M; L_0,L_1),
$$
 where the space $C^\infty(u^*(TM))$
is the space of vector fields along $u(\mathbb{R}\times [0,1])$, that is sections of
$u^*(TM)\to \mathbb{R}\times [0,1]$. Notice that the Cauchy-Riemann equation (\ref{e:hol})
defines a section, $u\mapsto \overline{\partial}_J(u)$, of this bundle. Moreover the 
moduli space $\hat{\mathcal{M}}_J(p,q,L_0,L_1)$ is precisely
the zero locus of this section. In order to show that the moduli space is a finite dimensional manifold, 
the section $\overline{\partial}_J$ must intersect transversally the zero-section.
Transversality is one of the problems in defining Lagrangian Floer homology,
it is  a delicate issue of the subject. 

The issue of  transversality of the section $\overline{\partial}_J$ is in fact relaxed, 
from the one stated above in the sense that the smooth condition on the map $u$ is relaxed.
The smooth condition is weakened to the Sobolev space $W^k_p(\mathbb{R}\times [0,1], M; L_0,L_1)$
for $k>p/2$ and $p>1$.
However, the new zero locus obtained in this setting coincides with the previous one due to
elliptic regularity; that is if $u\in W^k_p(\mathbb{R}\times [0,1], M; L_0,L_1)$ is such that
$\overline{\partial}_J(u)=0$, the $u$ is in fact smooth.

The main result in this direction is that there exists  a dense subset 
$\mathcal{J}_\textup{reg}(L_0,L_1)$
of $C^\infty([0,1],\mathcal{J})$   of $\omega$-compatible almost complex structures
such that for $J=\{J_t\} \in\mathcal{J}_\textup{reg}(L_0,L_1)$
and every $u\in \hat{\mathcal{M}}_J(p,q,L_0,L_1)$ the linearized operator  
$$
D(\overline{\partial}_J)_u: 
\{     \xi\in W^k_p(u^*(TM ))| \xi(s,0) \in L_0\textup{ and }  \xi(s,1) \in L_1\}
\to W^{k-1}_p(u^*(TM ))
$$
is a surjective Fredholm operator.
Furthermore the index of the operator $D(\overline{\partial}_J)_u$ is
the Maslov index of the map $u:\mathbb{R}\times [0,1]\to (M,\omega)$. 
Below we give the  definition of the Maslov index of $u$.
It then follows that the kernel of $D(\overline{\partial}_J)_u$ is finite dimensional
and  is identified with
the tangent space of $\hat{\mathcal{M}}_J(p,q,L_0,L_1)$ at $u$. For the details of these assertions
see \cite{floer-morse}.

As in the finite dimensional Morse theory case, the space of flow lines $
\hat{ \mathcal{M}}_J(p,q,L_0,L_1)$ admits an action of $\mathbb{R}$ on the $s$ coordinate. 
%To be specific,
%$\sigma.u(s,t):=u(s+\sigma,t)$ for $\sigma\in \mathbb{R}$. 
The quotient space by this action is denoted
 by ${\mathcal{M}}_J(p,q,L_0,L_1)$.
The space ${\mathcal{M}}_J(p,q,L_0,L_1)$ still needs to be taken  further apart; namely the homotopy
class of an element needs to be taken into consideration.
Let $\beta\in \pi_2(M,L_0\cup L_1 )$, and define ${\mathcal{M}}_J(p,q,L_0,L_1;\beta)$
as the elements $u\in {\mathcal{M}}_J(p,q,L_0,L_1)$ such that $[u]=\beta$.

In order to address the dimension of  ${\mathcal{M}}_J(p,q,L_0,L_1)$, 
for the moment consider only one Lagrangian submanifold $L$ of $(M,\omega)$. Then 
to each smooth map  $u:(D^2, \partial D^2)\to (M,L)$, 
one gets a trivial fibration $u^*(TM)\to D^2$  that is symplectic. Furthermore,
the fibration is trivial as symplectic bundles, $u^*(TM)\simeq \mathbb{R}^{2n}\times D^2$. 
Then when the fibration is restricted to $\partial D^2$, it 
defines a loop of Lagrangian subspaces of $(\mathbb{R}^{2n},\omega_0)$. 
That is, if $\Lambda(\mathbb{R}^{2n})$ represents the Grassmannian of Lagrangian subspaces
then the trivialized fibration induces a map $u_L:S^1\to \Lambda(\mathbb{R}^{2n})$.
The  {\em Maslov index}  $\mu_L(u)$,  of $u:(D^2, \partial D^2)\to (M,L)$ is defined to 
be the integer $(u_L)_*(1)\in \pi_1( \Lambda(\mathbb{R}^{2n}))\simeq \mathbb{Z}.$
The Maslov index of $u$ is well defined, it does not depend on the symplectic
trivialization; furthermore it only depends on the homotopy type
of $u$ relative to $L$.  Hence the Maslov index
induces a group morphism
$\mu_L:\pi_2(M,L)\to \mathbb{Z}$. 
The above concept extends to the case when two Lagrangian submanifolds 
 $L_0$ and $L_1$ are involved.
For the  definition of the Maslov index
see \cite{arnold-math}, \cite{ms} or \cite{robbin-salamon}. 

\begin{theorem}
\label{t:man}
Let $L_0$ and $L_1$ be compact Lagrangian submanifolds of $(M,\omega)$ that
%  $L_0\pitchfork L_1$.
intersect transversally.
Then there exists 
a dense subset $\mathcal{J}_\textup{reg}(L_0,L_1)$ 
in $C^\infty([0,1],\mathcal{J})$   of $\omega$-compatible almost complex structures, such that for 
$J=\{J_t\} \in \mathcal{J}_\textup{reg}(L_0,L_1)$,
$p$ and $q$ in $L_0\cap L_1$, and $\beta\in \pi_2(M,L_0\cup L_1 )$, the
space $\hat{\mathcal{M}}_J(p,q,L_0,L_1;\beta)$ is a smooth  manifold.
Moreover its dimension is given by the Maslov index
$\mu(\beta)$.
\end{theorem}

The proof of this result appears in \cite{floer-morse} and \cite{oh-floercoho}.
So far we imposed  conditions on $L_0, L_1$ and $(M,\omega)$, in order to have a more 
transparent
exposition of the subject. All the statements made so far hold for arbitrary closed symplectic 
manifolds and compact Lagrangian submanifolds, with the corresponding adaptations.
However the next results that we are going
to state does not hold  in  general. 
In fact it  is well understood that Lagrangian Floer 
homology can 
not be defined on  arbitrary symplectic manifolds for arbitrary Lagrangian submanifolds.

See \cite{fooo} and 
\cite{oh-floercoho} for more information on this peculiarity.

For instance one can impose the condition that the pair of  Lagrangians submanifolds of
$(M,\omega)$
must be monotone and that the  Maslov index of the Lagrangians has to be  greater than 2.
A Lagrangian submanifold $L\subset (M,\omega)$
is called {\em monotone} if there exists $\lambda>0$, such that
$\omega=\lambda\mu_L$.  One important fact that follows by considering monotone 
Lagrangian submanifold, 
is that the Maslov index of a non constant holomorphic disk with boundary in the 
Lagrangian is positive.
Under these conditions, Y.-G. Oh   \cite{oh-floercoho} defined 
Lagrangian Floer homology. 

\begin{comment}
BIEN ESTO QUEDA FUERA
A closed symplectic manifold $(M,\omega)$ is said to be monotone if there exists $\lambda>0$
such that $\omega=\lambda c_1$. Here $c_1$ is the first Chern class of  $(M,\omega)$ and
is computed with respect to any almost complex structure. 

Notice that if $L$ is a monotone Lagrangian submanifold, that is
there is $\lambda$ such that $\omega=\lambda\mu$, then $(M,\omega)$ must be
monotone and $\omega=(2\lambda)c_1$.
\end{comment}

There are less restrictive conditions for which Lagrangian Floer homology is well defined;
see \cite{fooo} and \cite{seidel-fukcat}. The advantage of the monotone assumptions is 
that the complex is just the $\mathbb{Z}_2$-vector space generated by the intersection points 
of the Lagrangian submanifolds; and the sum in the definition of the boundary operator (\ref{e:floerdiff}) 
is a finite sum. Basically the same picture as in the Morse theory case.

From now we are going to assume that the symplectic area of any 
2-disk with boundary in the
Lagrangian submanifold is zero. {\em Thus from now on we assume that
$(M,\omega)$ is a closed symplectic manifold and $L_0, L_1$ are closed Lagrangian submanifolds such that  
$[\omega]\cdot \pi_2(M,L_j)=0$ for $j=0,1$.} Under this conditions we will define the Floer
complex of the Lagrangian submanifolds $L_0$ and $L_1$.

Then 
under these hypothesis, when $\mu(p,q,\beta)=1$ the space 
${\mathcal{M}}_J(p,q,L_0,L_1;\beta)$ is compact and hence is a finite collection of points.
The reason why the space is compact is because 
 $[\omega]\cdot \pi_2(M,L_j)=0$ rules out the existence of holomorphic spheres and disks with
boundary in the Lagrangian submanifolds.
Since a convergent sequence of holomorphic strips, under Gromov's topology, converges
to a holomorphic strip with the possible union of holomorphic spheres and disks, the moduli
space is compact. See for example \cite{oh-floercoho}.
In some sence, this is most straightforward way scenario to define Lagrangian homology;
avoid holomorphic spheres and disks.

Now the advantage of using the field $\mathbb{Z}_2$, is that we don't have to worry about orientations
of the moduli spaces, that is assigning $(+)$ or $(-)$ to each component of 
 ${\mathcal{M}}_J(p,q,L_0,L_1;\beta)$
Recall that $J=\{J_t\}$ is given as in Theorem \ref{t:man}. 
Denote by $\#_{\mathbb{Z}_2 } {\mathcal{M}}_J(p,q,L_0,L_1;\beta)$ the number of points
of ${\mathcal{M}}_J(p,q,L_0,L_1;\beta)$
module 2.
Before defining the boundary operator as before, we note that the
solutions of (\ref{e:hol}) might determine an 
 infinite number of homotopy classes of
$\pi_2(M,L_0\cup L_1 )$. Thus we introduce the Novikov field over $\mathbb{Z}_2$
to give meaning to  the possible infinite number of homotopy classes of connecting  orbits.
The {\em Novikov field}
is defined as
$$
\Lambda:=  \left\{  	\sum_{j=0}^\infty  a_j T^{c_j} | a_j\in \mathbb{Z}_2, c_j\in \mathbb{R},
\lim_{j\to\infty} c_j=\infty\right \}.
$$
Let $\textup{CF}(L_0,L_1)$ be the free
$\Lambda$-module generated by the intersection points $L_0\cap L_1$, 
which are finitely many
since $L_0$ and $L_1$ are compact and intersect transversally. Then 
for $p\in L_0\cap L_1$  and $J=\{J_t\}$  as in Theorem \ref{t:man}, the boundary operator is defined as
\begin{eqnarray}
\label{e:floerdiff}
\partial_{J} (p):=  \sum_{
{\scriptsize
\begin{array}{c}
q\in L_0\cap L_1, \beta\in \pi_2(M,L_0\cup L_1) \\
\mu(p,q,\beta)=1
\end{array}}
}\hspace{ -1cm} \#_{\mathbb{Z}_2 } {\mathcal{M}}_J(p,q,L_0,L_1;\beta) \, \,  T^{\omega(\beta)} \, q.
\end{eqnarray}

\begin{theorem}
\label{t:diffloer}
Let $L_0$ and $L_1$ be closed Lagrangian submanifolds of  $(M,\omega)$ that intersect transversally
with $(M,\omega)$ also closed 
and  $\{J_t\}$ an  almost complex structure given by Theorem \ref{t:man}.
If $[\omega]\cdot \pi_2(M,L_j)=0$ for $j\in\{0,1\}$,
then the boundary operator $\partial_{J} :\textup{CF}(L_0,L_1) \to \textup{CF}(L_0,L_1)$ satisfies
\begin{eqnarray}
\label{e:dif}
\partial_{J} \circ\partial_{J}=0.
\end{eqnarray}
\end{theorem}

For the proof of this result see \cite{floer-morse}. For the  case of monotone 
Lagrangian submanifolds see \cite{oh-floercoho}.

Is important to know that the condition $[\omega]\cdot \pi_2(M,L_j)=0$ is fundamental in order
to have $
\partial_{J} \circ\partial_{J}=0.
$ 
In general Eq. (\ref{e:dif}) does not hold for arbitrary
 closed symplectic manifolds $(M,\omega)$ and compact  Lagrangian submanifolds $L_j$. 
In the context of Theorem \ref{t:diffloer},
the complex
$(\textup{CF}(L_0,L_1),\partial_{{J}})  $ is called the {\em Floer chain  complex}
of $(L_0,L_1)$. 
The {\em Lagrangian Floer Homology
of $(L_0,L_1)$} is defined to be the homology 
this complex,
$$
\textup{HF}(L_0,L_1,{J}) : = \frac{\textup{ Ker }\partial_{J}}{\textup{ Im }\partial_{J}}.
$$

The definition of the Floer differential is  similar to the one in Morse theory; in particular
in both cases we are counting the number of trajectories module 2. 
As pointed out 
before, it is possible to
use integer coefficients in the case of Morse homology. 
This is possible by fixing a coherent system
of orientation on each moduli space of gradient trajectories.
%we only have to fix an orientation on
%each unstable submanifolds for each non-degenerate fix point. 
However in the context of Lagrangian Floer
homology this is not the case;  the orientation issue is more involved in the Floer case.

Theorem \ref{t:diffloer} is the cornerstone of Lagrangian Floer homology. Here the hypothesis that 
the symplectic area of any 2-disk with boundary in an Lagrangian submanifold is crucial.
In fact there are known examples where Theorem  \ref{t:diffloer}  fails; \cite{oh-floercoho}, \cite[Ch. 2]{fooo}. 
The next result, Theorem \ref{t:varios},  is the philosophy of Lagrangian Floer homology.
Namely, is a blueprint to solve Arnol'd Conjecture that we will see below and  in  Section
\ref{s:app}.

\begin{theorem}
\label{t:varios}
Let $L_0,L_1$ and $(M,\omega)$  as in Theorem \ref{t:diffloer}.
Then 
\begin{itemize}
\item[(a)] $\quad\textup{HF}(L_0,L_1,J) $ is independent of $J=\{J_t\}\in \mathcal{J}_\textup{reg}(L_0,L_1)$,
\item[(b)] $\quad \textup{HF}(L_0,L_1) \simeq \textup{HF}(L_0,\phi(L_1))$ where $\phi$ is any Hamiltonian such that 
$L_0$ and  $\phi(L_1)$ intersect transversally, and 
\item[(c)] $\quad  \textup{HF}(L,L) \simeq H_*(L,\mathbb{Z}_2)\otimes_{\mathbb{Z}_2} \Lambda$. 
\end{itemize}
The isomorphisms in (b) and (c) are as $\Lambda$-modules.
\end{theorem}

In (c), $\textup{HF}(L,L)$ is understood in the sense of (b). That is, $\textup{HF}(L,L)$ is defined as
$\textup{HF}(L,\phi(L))$ where $\phi$ is any Hamiltonian diffeomorphisms such that $L$ and $\phi(L)$
 intersect transversally.

Theorems \ref{t:diffloer} and \ref{t:varios} reflect the idea of what is expected of 
Lagrangian Floer homology theory; in the sense that the theory must solve
Arnold's conjecture.
First, one requires the differential complex to be generated by the
intersection points of the Lagrangian submanifolds $L_0$ and $L_1$, and the differential 
operator  to square to zero.
Thus one has a homology theory, $\textup{HF}(L_0,L_1)$,  of the pair of Lagrangian submanifolds 
$(L_0,L_1)$. 
Finally one expects the  theory to 
satisfy (b) and (c) of Theorem \ref{t:varios}.  
With this at hand,  we have for any
 Lagrangian submanifold $L$ and
 any Hamiltonian diffeomorphism $\phi$, such that $L$ and $ \phi(L)$ intersect transversally, that 
\begin{eqnarray*}
\# (L\cap  \phi(L)) &\geq &\textup{dim}_\Lambda \textup{HF}(L, \phi(L))\\
&=&\textup{dim}_\Lambda \textup{HF}(L, L)\\
&= &\textup{dim}_\Lambda \Lambda\otimes_{\mathbb{Z}_2} H_*( L,\mathbb{Z}_2)\\
&= & \sum_{j=0}^{n} \textup{Rank }  H_j( L,\mathbb{Z}_2).\\
\end{eqnarray*}
where the Lagrangian $L$ has dimension $n$.

\begin{theorem}[A. Floer \cite{floer-morse}]
Let $(M,\omega)$ be a closed symplectic manifold, $L$ a Lagrangian submanifold
and $\phi$ a Hamiltonian diffeomorphisms such that $L$ and $\phi(L)$ intersect
transversally. Further, assume that $[\omega]\cdot \pi_2(M,L)=0$, then
$$
\# (L\cap  \phi(L))\geq \sum_{j=0}^{n} \textup{Rank }  H_j( L,\mathbb{Z}_2).
$$
\end{theorem}

This result was generalized by K. Fukaya, Y.-G. Oh, H. Otha and K. Ono in 
\cite{fooo}  using the technique of Kuranishi structures on the moduli
space. The hypothesis $[\omega]\cdot \pi_2(M,L)=0$ is replaced by the
requirement that the map $H_*(L,\mathbb{Q})\to H_*(M,\mathbb{Q}) $
must be injective. Notice that this new hypothesis can not be relaxed.
As mentioned in Section \ref{s:symandlag}, any symplectic manifold admits
arbitrary small Lagrangian tori that are displaceable.

Now in the context of Arnol'd Conjecture,  let $(M,\omega)$ be a closed symplectic 
manifold and  $\psi$ be a Hamiltonian diffeomorphism
with non-degenerate fixed points. Then we know that  $\textup{graph}(\psi)$ is a Lagrangian 
submanifold of the symplectic manifold $(M\times M, \pi_1^*(\omega)-\pi_2^*(\omega))$.
Since the fixed points of $\psi$ are non-degenerate,  the intersection of 
$\textup{graph}(\psi)$ with 
 the diagonal $\Delta$ is transversal. Now assume that 
$(M\times M, \pi_1^*(\omega)-\pi_2^*(\omega))$
and $\textup{graph}(\psi)$ satisfy the hypotheses about the zero symplectic area of any 2-disk with boundary
in the Lagrangian, then the above computation implies that 
\begin{eqnarray*}
\# (\Delta \cap  \textup{graph}(\psi)) &=& 
\# (\Delta \cap  (1\times\psi)(\Delta)) \\
&\geq &\!
\sum_{j=0}^{2n} \textup{Rank }  H_j( \Delta,\mathbb{Z}_2)\\
&=&\!\sum_{j=0}^{2n} \textup{Rank }  H_j( M,\mathbb{Z}_2).
\end{eqnarray*}
That is $|\textup{Fix}(\psi)| \geq\!\sum_{j=0}^{2n} \textup{Rank }  H_j( M,\mathbb{Z}_2)$.

\begin{corollary}[A. Floer \cite{floer-morse}]
Let $(M,\omega)$ be a closed symplectic manifold and  $\phi$ a Hamiltonian 
diffeomorphism with non-degenerate critical points. If $\pi_2(M)=0$, then
$$
|\textup{Fix}(\phi)| \geq\!\sum_{j=0}^{2n} 
\textup{Rank }  H_j( M,\mathbb{Z}_2).
$$
\end{corollary}

\begin{comment}
$$
\# (\Delta \cap  \textup{graph}(\psi)) = 
\# (\Delta \cap  (1\times\psi)(\Delta)) \geq \!
\sum_{j=0}^{2n} \textup{Rank }  H\!\!\!_j\!( \Delta,\mathbb{Z}_2)
=\!\sum_{j=0}^{2n} \textup{Rank }  H\!\!\!_j\!( M,\mathbb{Z}_2).
$$
\end{comment}

\medskip
These notes are far from been an introduction to the subject of Lagrangian Floer homology.
There are many important issues that we did not mentioned at all. For example the regularity
of holomorphic strip, and 
 the transversality issue  to assure that the space
of trajectories $\hat{\mathcal{M}}(p,q;\beta)$ is a smooth manifold. Another important fact
that we left out was the compactness issue of $\hat{\mathcal{M}}(p,q;\beta)$, among other things. 
See \cite{audin-morseth},
\cite{fooo} and \cite{seidel-fukcat}.  
Finally we would like to mention that  Lagrangian Floer homology admits a product structure
$$
\textup{HF}(L_0,L_1)\otimes \textup{HF}(L_1,L_2)\to \textup{HF}(L_0,L_2).
$$
That is, instead of considering only two Lagrangian submanifolds one considers three 
Lagrangians.
Instead of considering holomorphic strips, one considers holomorphic triangles, that is map
from the closed unit disk minus three point of the boundary; and each component of the boundary is
mapped  to a different Lagrangian submanifold. 
Of course there is nothing special of taking three Lagrangian
submanifolds, one can consider any number of Lagrangian submanifolds and the corresponding holomorphic
polygons. 
In formal terms, is said that Lagrangian Floer homology admits an $A_\infty$-structure.
For further reading of  this structure see \cite{fooo} and \cite{seidel-fukcat}. 
Also \cite{auroux-a} for an introduction into the subject.

%%%%%%%%%%%%%%%%%%%%%%%%%%%%%%%%%%%%%%%%%%%%%%%%%%%%%%%%%%
%%%%%%%%%%%%%%%%%%%%%%%%%%%%%%%%%%%%%%%%%%%%%%%%%%%%%%%%%%
\section{Computation of $\textup{HF}(L,L)$}
%%%%%%%%%%%%%%%%%%%%%%%%%%%%%%%%%%%%%%%%%%%%%%%%%%%%%%%%%%
%%%%%%%%%%%%%%%%%%%%%%%%%%%%%%%%%%%%%%%%%%%%%%%%%%%%%%%%%%

In this section we  give the outline that shows that
$\textup{HF}(L,L)=\textup{HF}(L,\phi(L))$ is isomorphic to 
$H_*(L,\mathbb{Z}_2)\otimes_{\mathbb{Z}_2} \Lambda$ under the assumption that 
$[\omega]\cdot \pi_2(M,L)=0$. As pointed out in the previous section, this 
isomorphism is fundamental
in the proof of Arnol'd Conjecture. 
The idea behind the proof of
$\textup{HF}(L,\phi(L)) \simeq H_*(L,\mathbb{Z}_2)\otimes_{\mathbb{Z}_2} \Lambda$
is to prove it in the case when the symplectic manifold is $(T^*L, -d\lambda_{\textup{can}})$
and the Lagrangian submanifold is the zero section $L_0$. That is 
$\textup{HF}(L_0,\phi(L_0)) \simeq H_*(L_0,\mathbb{Z}_2)\otimes_{\mathbb{Z}_2} \Lambda$
for some particular Hamiltonian $\phi:T^*L\to T^*L$.
Though in this case the symplectic manifold is non compact, 
 together with the hypothesis
 $[\omega_\textup{can}]\cdot \pi_2(T^*L,L_0)=0$ we assume that 
  Theorem \ref{t:diffloer} still hols  in this case.

The proof  of the general case, that is for an arbitrary compact symplectic manifold
 $(M,\omega)$ and a Lagrangian submanifold $L$  relies on Weinstein's Lagrangian tubular neighborhood
theorem. That is,
consider a Hamiltonian $\phi: M\to M$ such that $\phi(L)$ lies in a Weinstein's tubular neighborhood
of $L$. Then under the symplectic diffeomorphism given by Weinstein theorem, 
the known computations of $(T^*L, -d\lambda_{\textup{can}})$, that is holomorphic
disks and gradient flow lines in $L$, are carried to
 $(M,\omega)$ and $L$ to obtain the isomorphism between $\textup{HF}(L,\phi(L))$
 and $H_*(L,\mathbb{Z}_2)\otimes_{\mathbb{Z}_2} \Lambda$.

Now we show the Riemannian structure and almost complex structure on the cotangent
bundle that will be relevant in this section in the computation of Lagrangian Floer homology.
Hence given a Riemannian metric $g$ on $L$
%%%%%%%%% así es ????
it induces an almost complex structure ${\bf J}$ on $T^*L$ that is compatible with
$\omega_\textup{can}$. Such that at the zero section $(p,0)\in T^*L $,
it maps vertical vectors $T^*_pL\subset T_{(p,0)}(T^*L)$ to horizontal vectors
$T_pL\subset T_{(p,0)}(T^*L)$. Set ${\bf g}$ to be  the Riemannian structure induced by 
${\bf J}$  and $\omega_\textup{can}$, thus ${\bf g}(\cdot,\cdot) = \omega_\textup{can}(\cdot, {\bf J} \cdot  )$.
Notice that such almost complex structure  ${\bf J}$ is not unique.

Now lets see the relation between Morse homology on $L$ and Lagrangian Floer homology
of the zero section $L_0$
in $(T^*L, \omega_\textup{can})$.
To that end let $g$ be a Riemannian structure  on $L$
and  $f:L\to\mathbb{R}$ be a Morse-Smale function.
Then the graph of the 1-form $df$,
$$
L_1:=\textup{graph}(df)=\{(p, df_p)| p\in L   \},
$$
is a Lagrangian submanifold in  $(T^*L, \omega_{\textup{can}})$ since is 
the graph of a closed 1-form. 
%%%% aquí voy
%Denote by $L_0$ be the zero section in $(T^*L, \omega_{\textup{can}})$. 
Therefore 
 $L_1$ intersects $L_0$  precisely at the the critical points of $f$; 
that is at points $(p,0)$ for $p\in \crit{f}$. 
 Furthermore 
 the intersection is transversal and consists of finitely many points. 

As before 
let  $\pi: T^*L\to L$ be the projection map,
set $H:=-f\circ \pi$ and  $X_H$ the Hamiltonian vector field of $H$ on 
$(T^*L, \omega_{\textup{can}})$, thus
$$
\omega_{\textup{can}}(X_H,\cdot )=dH.
$$ 
Notice that if $(x_1,\ldots, x_n )$ are the local coordinates of $L$, and 
$(x_1,\ldots, x_n,y_1,\ldots ,y_n )$ the corresponding local coordinantes
of $T^*L$ as in Example \ref{exa:cotang}, then the Hamiltonian vector
field $X_H$ takes the form,
$$
X_H:= \frac{\partial f}{\partial x_1}\frac{\partial }{\partial y_1} + \cdots +
\frac{\partial f}{\partial x_n} \frac{\partial }{\partial y_n}.
$$
This expression justifies the minus sign in the definition of the Hamiltonian function
$H$.
Finally let  $\phi_t:T^*L\to T^*L$ be the  path of Hamiltonian  diffeomorphisms
induced by the Hamiltonian function $H$.

In local coordinates $(x_1,\ldots, x_n,y_1,\ldots ,y_n )$, the path 
$\{\phi_t\}$ of Hamiltonian  diffeomorphisms takes the form
$$
 \phi_t(x_1,\ldots, x_n,y_1,\ldots ,y_n )=
 \left(x_1,\ldots, x_n, t \frac{\partial f}{\partial x_1}(x)+y_1   ,\ldots ,t \frac{\partial f}{\partial x_n}(x) 
 +y_n\right ).
$$
Therefore $\phi_1(p,0)=(p,df_p)$.
\begin{proposition}
Let $\phi_t:T^*L\to T^*L$, $L_0$ and $L_1$ as above. Then the time-1 map
is such that  $\phi_1(L_0)=L_1.$ 
\end{proposition}

However to describe $X_H$ in global terms we must use the Riemannian structure $g$ on $L$
and the remark made above. Thus let ${\bf J}$ on $(T^*L, \omega_{\textup{can}})$ as before
subject to the condition that on the zero section 
 ${\bf J} \,\textup{grad}(f)= df\in T^*L\subset T(T^*L)$ where
$\textup{grad}(f)\in TL\subset T(T^*L)$. Furthermore  on the zero section we have that
$$
\textup{{\bf grad}} (H)= \textup{grad}(f).
$$
Since ${\bf J}$ is compatible with  $\omega_{\textup{can}}$
it follows from $\omega_{\textup{can}}(X_H,\cdot )=dH$ that 
$\textup{{\bf grad}} (H)= -{\bf J} \, X_H$ and $X_H={\bf J} \,\textup{grad}(f)$.

In the context of Lagrangian Floer homology we care about the intersection points
of $L_0$ with $\phi(L_0)$, that in this particular case  are in a  natural bijection
with the critical points
of $f$. Therefore there is a natural $\Lambda$-linear bijection
$$
I:\textup{CF}(L_0,\phi(L_0)) \to \textup{Crit}(f)\otimes_{\mathbb{Z}_2} \Lambda
$$
which on generators takes the form $I(p,0)=p$. This map will give rise to the
isomorphism between 
$\textup{HF}(L_0,\phi(L_0))$ and $ H_*(L,\mathbb{Z}_2)\otimes_{\mathbb{Z}_2} \Lambda$.
It remains to look at the differentials maps  in each case; and henceforth holomorphic 
disks and gradient flow lines.

However in order to  relate the flow lines  in $L$  of the gradient of $f$  of Morse index 1 that 
come into play in the 
Morse differential, with the holomorphic disks in $T^*L$ with  boundary  in $L_0$ and $ $ and of
Maslov index 1 that come into play in the Floer differential, and additional condition on
$f$ is imposed. %To that end we must restrict to a particular class of functions $f$. 
Namely $f$ must be $C^2$-small with respect to $g$.
Recall that we assume that $f:N\to\mathbb{R}$ satisfies the Morse-Smale condition.  
Then according to A. Floer \cite{floer-witten},  the family of  almost complex structures 
$
J_t:=- \phi_{t,*} \circ {\bf J} \circ  (\phi_{t,*})^{-1}
$
for $0\leq t\leq 1$,  is regular. That is  $J=\{J_t\} \in \mathcal{J}_\textup{reg}(L_0,L_1)$.
Since we assume that $[\omega_\textup{can}]\cdot\pi_2(T^*L,L_0)=0$, then the Lagrangian Floer
homology of $(L_0,L_1)$ can be computed using the moduli spaces $ {\mathcal{M}}_J(p,q,L_0,L_1)$,
where $p$ and $q$ are intersection points.
Recall that for $p$ and $q$ critical points of $f:L\to\mathbb{R}$, 
$ \mathcal{M}(f;p,q)$ is the moduli space of flow lines of $-\textup{grad}(f)$ connecting
$p$ to $q$.

Consider the map %Furthermore there is a bijection
$$
\Psi:     \mathcal{M}(f;p,q)  \to  {\mathcal{M}}_J(p,q,L_0,L_1)
$$
given by $\Psi(u)(s,t):= \phi_t(u(s)).$ 
Next we show that $\Psi$ is well defined, that is that $\phi_t(u(s))$ is $J$-holomorphic
and satisfies the boundary conditions. For
note that $\Psi(u)(s,0)=\phi_0(u(s))$ lies in $L_0$
and $\Psi(u)(s,1)=\phi_1(u(s))$ lies in $L_1$. It only remains to show that $\Psi(u)$ is
$J$-holomorphic. To that end notice that  
%If $u_0$ is in $ \mathcal{M}(f;p,q)$, in particular we have that 
$$
\frac{d u}{dt}(t)=-\grad{f}(u(t)), 
$$
Therefore
\begin{eqnarray*}
\frac{\partial}{\partial s }\Psi(u)(s,t) &=&\frac{\partial}{\partial s }  \phi_t(u(s))
=\phi_{t,*}(u(s))  \left( \frac{d}{d s }  u(s)\right) \\
&=&-\phi_{t,*} (u(s))\left( \grad{f}({u(t))} \right) .
\end{eqnarray*}
On the other hand, since $X_H={\bf J} \,\textup{grad}(f)$  and $-J_t\circ  \phi_{t,*}= \phi_{t,*} \circ {\bf J} 
$ it follows that  
\begin{eqnarray*}
\frac{\partial}{\partial t }\Psi(u)(s,t) &=&\frac{\partial}{\partial t }  \phi_t(u(s))
= X_H(   \phi_t( u(s))  )   \\ %%   (\phi_t(u_0(s))) \\
&=&\phi_{t,*}(u(s))  (X_H (u(s))  )\\
&=& \phi_{t,*}(u(s))   ({\bf J} \;  \textup{grad}(f)(u(s)))\\
&=&- J_t \phi_{t,*}(u(s))   (\textup{grad}(f)(u(s)))  %%% J_{t, \phi_t(u_0(s))}      \ \   \ \ 
\end{eqnarray*}
Since  $J_t^2=-1$, we get that 
\begin{eqnarray}
%\label{e:hol}
\frac{\partial}{\partial s }\Psi(u)(s,t)  +J_t
\frac{\partial}{\partial t }\Psi(u)(s,t) 
=0.
\end{eqnarray}
 That is if $u$ is a gradient flow line we have that $\Psi(u)$ satisfies the Cauchy-Riemann
equation with respecto to $J=\{J_t\}$. That is, $\Psi:     \mathcal{M}(f;p,q)  
\to  {\mathcal{M}}_J(p,q,L_0,L_1)$ is well-defined and according to A. Floer \cite{floer-witten} is a bijection
for $f$ small enough.

Recall that 
$$
I:\textup{CF}(L_0,\phi(L_0))  \to \textup{Crit}(f)\otimes_{\mathbb{Z}_2} \Lambda
$$
given by $I(p,0)=p$ is also a bijection.  Furthermore, the fact that $\Psi$
is a bijection implies that $I$ is a chain map, $I\circ \partial_J =\partial\circ I$. 
(Here the Morse differential on $\textup{Crit}(f)$ is extended $\Lambda$-linearly to 
$\textup{Crit}(f)\otimes_{\mathbb{Z}_2} \Lambda$)
Moreover
the induced map $I: \textup{HF}(L_0,\phi(L_0))\to \textup{MH}_*(L;f,g) \otimes_{\mathbb{Z}_2} \Lambda  $
is an isomorphism of $\Lambda$-modules. Finally, since  $\textup{MH}_*(L;f,g)\simeq H_*(L,\mathbb{Z}_2) $
we have that 
\begin{eqnarray}
\label{e:HFiso}
\textup{HF}(L_0,\phi(L_0))
\simeq H_*(L,\mathbb{Z}_2)\otimes_{\mathbb{Z}_2} \Lambda .
\end{eqnarray}
in the case of the symplectic conagent bundle $(T^*L, \omega_{\textup{can}})$
and the zero section $L_0$ Lagrangian submanifold.

The proof of the above statement in the case of an arbitrary compact
 symplectic manifold $(M,\omega)$ and  a compact  
Lagrangian submanifold $L$,
under the assumption that $[\omega]\cdot \pi_2(M,L)=0$,  is based on the 
on the contagent bundle case. For, consider $\phi:M\to M$ a Hamiltonian diffeomorphisms
such that $L$ and $\phi(L)$ intersect transversally and  small enough so that
$\phi(L)$ lies in tubular neighborhood of $L$. Then by 
Weinstein's Lagrangian neighborhood
theorem, this small neighborhood of $L$ is symplectomorphic to 
a neighborhood to the zero section of $(T^*L,\omega_{can})$.
Hence the all the relevant information
in the computation of $\textup{HF}(L,\phi(L))$    lie in the tubular neighborhood of $L$. 
Thus the isomorphism (\ref{e:HFiso}) also holds in this case.

%%%%%%%%%%%%%%%%%%%%%%%%%%%%%%%%%%%%%%%%%%%%%%%%%%%%%%%%%%
%%%%%%%%%%%%%%%%%%%%%%%%%%%%%%%%%%%%%%%%%%%%%%%%%%%%%%%%%%

\section{Applications}
\label{s:app}
%%%%%%%%%%%%%%%%%%%%%%%%%%%%%%%%%%%%%%%%%%%%%%%%%%%%%%%%%%
%%%%%%%%%%%%%%%%%%%%%%%%%%%%%%%%%%%%%%%%%%%%%%%%%%%%%%%%%%

In principle Lagrangian Floer homology, was meant to solve Arnol'd conjecture. 
Nowadays it is important on its own. In this section we will briefly 
explain a few examples.

The 2-sphere example revisited. Consider $(S^2,\omega)$ and $L$ any embedded 
circle. Clearly this example does
not satisfy the condition $[\omega]\cdot \pi_2(S^2,L)=0$; nevertheless for 
some particular Lagrangian submanifolds it lies in the monotone case. For instance if
$L$ is an equator then is a monotone Lagrangian and the Lagrangian Floer homology applies.
Consider $\phi$ is a rotation of $(S^2,\omega)$ such that 
$L$ and $\phi(L)$ are transversal, then the standard complex structure $J$ %on $(S^2,\omega)$
is regular. In this case $\textup{CF}(L,\psi(L))$ has two generators and 
$\partial_J$ is the zero map. Therefore $\textup{HF}(L,\phi(L))$ has rank two. 
As we explain above, this means that the equator is a non displaceable 
Lagrangian of $(S^2,\omega)$.

More generally, consider the Lagrangian submanifold
$\mathbb{R}P^n$  in 
  $(\mathbb{C}P^n,\omega_{FS})$ for $n\geq 1$. The Lagrangian  
$\mathbb{R}P^n$ is a monotone. On   $(\mathbb{C}P^n,\omega_{FS})$ there
is a canonical $SU(n+1)$ Hamiltonian action, that is induced from the 
standard  linear action on the euclidean space $\mathbb{C}^{n+1}.$
Restrict the action to the maximal torus of $SU(n+1)$. Then a vector
on the Lie algebra of the maximal torus induces, by the exponential map,
a Hamiltonian diffeomorphism $\phi$ of  $(\mathbb{C}P^n,\omega_{FS})$.
Moreover, is such that  $\phi(\mathbb{R}P^n)$  meets $\mathbb{R}P^n$
transversally. In this example the standard complex structure $J$,
$\mathbb{R}P^n$ and  $\phi(\mathbb{R}P^n)$ satisfy the 
transversality condition, that is the moduli spaces  
 ${\mathcal{M}}_J(p,q;\beta)$ are smooth manifolds for $p,q\in 
\mathbb{R}P^n\cap\phi(\mathbb{R}P^n)$. Moreover 
the analog of Theorem \ref{t:diffloer} in the monotone case
also holds. 
Hence the Lagrangian Floer homology of
$( \mathbb{R}P^n,\phi(\mathbb{R}P^n))$ is well defined, further
the differential $\partial_J$ is the zero map and 
$$
\textup{HF}( \mathbb{R}P^n,\phi(\mathbb{R}P^n))
$$
has rank $n+1$ as a  $\Lambda$-module. In particular, $\mathbb{R}P^n$ is non displaceable. 
This result is due to Y.-G. Oh \cite{oh-floerii}, in the setting of 
Lagrangian Floer
homology for monotone Lagrangians.

Another important example concerns symplectic toric manifolds. In this case
if $(M,\omega)$ is a toric manifold with moment map $\Phi:M\to \Delta$, then 
for $x$ in the interior of the polytope $\Delta$, $L_x=\Phi^{-1}(x)$ is a Lagrangian
$n$-torus. In this case there is a dichotomy,
$$
\textup{HF}(L_x,L_x)=0\;\;\;\; \textup{ or } \;\;\;\;\textup{HF}(L_x,L_x)=
H_*(\mathbb{T}^n;\mathbb{Z}_2) \otimes \Lambda.
$$
In fact a stronger result is true, those Lagrangian torus $L_x$ such that
$\textup{HF}(L_x,L_x)=0$ are in fact displaceable.

Example \ref{exa:sphere} of the 2-sphere with a Morse function $f$, is an example 
of a symplectic toric manifold. In this case the Morse function is in fact
a  moment map, recall that $f:S^2\to [0,1]$
is given by $f(x,y,z)=z$. Then for $z\in [0,1]$ not equal to zero
$L_z$ is a circle that lies entirely in the north or south hemisphere; hence $L_z$
is displaceable. And for the equator $L_0$,  we have
$$
\textup{HF}(L_0,L_0)=
H_*(S^1;\mathbb{Z}_2) \otimes \Lambda,
$$
as mentioned at the beginning of this section. For further details
on Lagrangian Floer homology on symplectic toric manifolds 
see \cite{fooo-torici}.

\begin{comment}
%%%%%%%%%%%%%%%%%%%%%%%%%%%%%%%%%%%%%%%%%%%%
%%%%%%%%%%%%%%%%%%%%%%%%%%%%%%%%%%%%%%%%%%%%%%%%%%%%%%%%%%
\section{H(L,L)}
%%%%%%%%%%%%%%%%%%%%%%%%%%%%%%%%%%%%%%%%%%%%
%%%%%%%%%%%%%%%%%%%%%%%%%%%%%%%%%%%%%%%%%%%%%%%%%%%%%%%%%%

The way to compute $HF(L,L)$, consists in taking a Hamiltonian diffeomorphisms $\psi$
such that $\psi(L)$ and $L$ intersect transversally. Then   $HF(L,L)=HF(L,\psi(L))$.

Let $L$ be a Lagrangian submanifold of $(M,\omega)$. Then by Weinstein's Lagrangian tubular neighborhood theorem
there are  neighborhoods $U$ of $L\subset (M,\omega)$ and $V$ of the zero section $L_0$ in $(T^*L,
\omega_{\textup{can}})$ that are  diffeomorphic  by a map $\phi$ such that $\phi^*(\omega_{\textup{can}}
)=\omega$. Now consider a smooth function $f:L\to\mathbb{R}$, then the graph of $df$ is a Lagrangian submanifold
$L_1$ in   $(T^*L,\omega_{\textup{can}})$. Moreover $L_0$ is Hamiltonian isotopic to $L_1$; for consider the
Hamiltonian function $H:=(\epsilon f)\circ \pi$. Furthermore, the intersection points of $L_0$ with $L_1$ correspond
to the critical points of $f$. Finally, choose $f$ so that $L_1\subset V$ and $L_0$ intersects transversally with $L_1$.

Consider $f:L\to\mathbb{R}$ be a smooth $C^2$-small function. And consider the 
time-1 Hamiltonian $\phi$ induced by $H:=f\circ \pi$. Then the Lagrangian 
$\phi(L_0)$ intersects $L_0$ transversally, and the intersection points correspond
to the critical points of $f$.
\end{comment}

\bibliographystyle{acm}
\bibliography{/Users/Andres/Dropbox/Documentostex/Ref}    %%%lap

\begin{thebibliography}{10}

\bibitem{arnold-sur}
{\sc Arnol'd, V.}
\newblock Sur une propri\'et\'e topologique des applications globalement
  canoniques de la m\'ecanique classique.
\newblock {\em C. R. Acad. Sci. Paris 261\/} (1965), 3719--3722.

\bibitem{arnold-math}
{\sc Arnol'd, V.~I.}
\newblock {\em Mathematical methods of classical mechanics}, vol.~60 of {\em
  Graduate Texts in Mathematics}.
\newblock Springer-Verlag, New York, 1993.
\newblock Translated from the 1974 Russian original by K. Vogtmann and A.
  Weinstein, Corrected reprint of the second (1989) edition.

\bibitem{audin-morseth}
{\sc Audin, M., and Damian, M.}
\newblock {\em Morse theory and {F}loer homology}.
\newblock Universitext. Springer, London; EDP Sciences, Les Ulis, 2014.
\newblock Translated from the 2010 French original by Reinie Ern{\'e}.

\bibitem{auroux-a}
{\sc Auroux, D.}
\newblock A beginner's introduction to {F}ukaya categories.
\newblock In {\em Contact and symplectic topology}, vol.~26 of {\em Bolyai Soc.
  Math. Stud.} J\'anos Bolyai Math. Soc., Budapest, 2014, pp.~85--136.

\bibitem{banyaga-surla}
{\sc Banyaga, A.}
\newblock Sur la structure du groupe des diff\'eomorphismes qui pr\'eservent
  une forme symplectique.
\newblock {\em Comment. Math. Helv. 53}, 2 (1978), 174--227.

\bibitem{bott-morse}
{\sc Bott, R.}
\newblock Morse theory indomitable.
\newblock {\em Inst. Hautes \'Etudes Sci. Publ. Math.}, 68 (1988), 99--114
  (1989).

\bibitem{cannas-lectures}
{\sc Cannas~da Silva, A.}
\newblock {\em Lectures on symplectic geometry}, vol.~1764 of {\em Lecture
  Notes in Mathematics}.
\newblock Springer-Verlag, Berlin, 2001.

\bibitem{conley-zehnder-thebir}
{\sc Conley, C.~C., and Zehnder, E.}
\newblock The {B}irkhoff-{L}ewis fixed point theorem and a conjecture of {V}.
  {I}. {A}rnol'd.
\newblock {\em Invent. Math. 73}, 1 (1983), 33--49.

\bibitem{eliashberg-atheorem}
{\sc Eliashberg, Y.~M.}
\newblock A theorem on the structure of wave fronts and its application in
  symplectic topology.
\newblock {\em Funktsional. Anal. i Prilozhen. 21}, 3 (1987), 65--72, 96.

\bibitem{floer-morse}
{\sc Floer, A.}
\newblock Morse theory for {L}agrangian intersections.
\newblock {\em J. Differential Geom. 28}, 3 (1988), 513--547.

\bibitem{floer-witten}
{\sc Floer, A.}
\newblock Witten's complex and infinite-dimensional {M}orse theory.
\newblock {\em J. Differential Geom. 30}, 1 (1989), 207--221.

\bibitem{fortune-weinstein-asymp}
{\sc Fortune, B., and Weinstein, A.}
\newblock A symplectic fixed point theorem for complex projective spaces.
\newblock {\em Bull. Amer. Math. Soc. (N.S.) 12}, 1 (1985), 128--130.

\bibitem{fukaya-morsehomotopy}
{\sc Fukaya, K.}
\newblock Morse homotopy, {$A^\infty$}-category, and {F}loer homologies.
\newblock In {\em Proceedings of {GARC} {W}orkshop on {G}eometry and {T}opology
  '93 ({S}eoul, 1993)\/} (1993), vol.~18 of {\em Lecture Notes Ser.}, Seoul
  Nat. Univ., Seoul, pp.~1--102.

\bibitem{fooo}
{\sc Fukaya, K., Oh, Y.-G., Ohta, H., and Ono, K.}
\newblock {\em Lagrangian intersection {F}loer theory: anomaly and obstruction.
  {P}art {I}}, vol.~46 of {\em AMS/IP Studies in Advanced Mathematics}.
\newblock American Mathematical Society, Providence, RI; International Press,
  Somerville, MA, 2009.

\bibitem{fooo-torici}
{\sc Fukaya, K., Oh, Y.-G., Ohta, H., and Ono, K.}
\newblock Lagrangian {F}loer theory on compact toric manifolds. {I}.
\newblock {\em Duke Math. J. 151}, 1 (2010), 23--174.

\bibitem{fukaya-ono-arnold}
{\sc Fukaya, K., and Ono, K.}
\newblock Arnold conjecture and {G}romov-{W}itten invariant.
\newblock {\em Topology 38}, 5 (1999), 933--1048.

\bibitem{gompf-symplectically}
{\sc Gompf, R.~E.}
\newblock Symplectically aspherical manifolds with nontrivial {$\pi_2$}.
\newblock {\em Math. Res. Lett. 5}, 5 (1998), 599--603.

\bibitem{gromov-psudo}
{\sc Gromov, M.}
\newblock Pseudoholomorphic curves in symplectic manifolds.
\newblock {\em Invent. Math. 82}, 2 (1985), 307--347.

\bibitem{hofer-salamon-floerhomo}
{\sc Hofer, H., and Salamon, D.~A.}
\newblock Floer homology and {N}ovikov rings.
\newblock In {\em The {F}loer memorial volume}, vol.~133 of {\em Progr. Math.}
  Birkh\"auser, Basel, 1995, pp.~483--524.

\bibitem{liu-tian-floerhomo}
{\sc Liu, G., and Tian, G.}
\newblock Floer homology and {A}rnold conjecture.
\newblock {\em J. Differential Geom. 49}, 1 (1998), 1--74.

\bibitem{matsumoto-morse}
{\sc Matsumoto, Y.}
\newblock {\em An introduction to {M}orse theory}, vol.~208 of {\em
  Translations of Mathematical Monographs}.
\newblock American Mathematical Society, Providence, RI, 2002.
\newblock Translated from the 1997 Japanese original by Kiki Hudson and
  Masahico Saito, Iwanami Series in Modern Mathematics.

\bibitem{mcduff-notes}
{\sc McDuff, D.}
\newblock Notes on kuranishi atlases.
\newblock {\em {\textit arXiv:1411.4306}\/}.

\bibitem{ms}
{\sc McDuff, D., and Salamon, D.}
\newblock {\em Introduction to symplectic topology}, second~ed.
\newblock Oxford Mathematical Monographs. The Clarendon Press, Oxford
  University Press, New York, 1998.

\bibitem{msjholo}
{\sc McDuff, D., and Salamon, D.}
\newblock {\em {$J$}-holomorphic curves and symplectic topology}, second~ed.,
  vol.~52 of {\em A. M. S. Colloquium Publications}.
\newblock American Mathematical Society, Providence, RI, 2012.

\bibitem{mw-2}
{\sc McDuff, D., and Wehrheim, K.}
\newblock The fundamental class of smooth kuranishi atlases with trivial
  isotropy.
\newblock {\em {\textit arXiv:1508.01560}\/}.

\bibitem{mw-1}
{\sc McDuff, D., and Wehrheim, K.}
\newblock Smooth kuranishi atlases with isotropy.
\newblock {\em {\textit arXiv:1508.01556 }\/}.

\bibitem{mw-3}
{\sc McDuff, D., and Wehrheim, K.}
\newblock The topology of kuranishi atlases.
\newblock {\em {\textit arXiv:1508.01844}\/}.

\bibitem{milnor-morse}
{\sc Milnor, J.}
\newblock {\em Morse theory}.
\newblock Based on lecture notes by M. Spivak and R. Wells. Annals of
  Mathematics Studies, No. 51. Princeton University Press, Princeton, N.J.,
  1963.

\bibitem{liviu-morse}
{\sc Nicolaescu, L.~I.}
\newblock {\em An invitation to {M}orse theory}.
\newblock Universitext. Springer, New York, 2007.

\bibitem{oh-floercoho}
{\sc Oh, Y.-G.}
\newblock Floer cohomology of {L}agrangian intersections and pseudo-holomorphic
  disks. {I}.
\newblock {\em Comm. Pure Appl. Math. 46}, 7 (1993), 949--993.

\bibitem{oh-floerii}
{\sc Oh, Y.-G.}
\newblock Floer cohomology of {L}agrangian intersections and pseudo-holomorphic
  disks. {II}. {$({\bf C}{\rm P}^n,{\bf R}{\rm P}^n)$}.
\newblock {\em Comm. Pure Appl. Math. 46}, 7 (1993), 995--1012.

\bibitem{oh-symplectic1}
{\sc Oh, Y.-G.}
\newblock {\em Symplectic topology and {F}loer homology. {V}ol. 1}, vol.~28 of
  {\em New Mathematical Monographs}.
\newblock Cambridge University Press, Cambridge, 2015.
\newblock Symplectic geometry and pseudoholomorphic curves.

\bibitem{oh-symplectic2}
{\sc Oh, Y.-G.}
\newblock {\em Symplectic topology and {F}loer homology. {V}ol. 2}, vol.~29 of
  {\em New Mathematical Monographs}.
\newblock Cambridge University Press, Cambridge, 2015.
\newblock Floer homology and its applications.

\bibitem{ono-onthe}
{\sc Ono, K.}
\newblock On the {A}rnol'd conjecture for weakly monotone symplectic manifolds.
\newblock {\em Invent. Math. 119}, 3 (1995), 519--537.

\bibitem{pardon-analgebraic}
{\sc Pardon, J.}
\newblock An algebraic approach to virtual fundamental cycles on moduli spaces
  of pseudo-holomorphic curves.
\newblock {\em Geom. Topol. 20}, 2 (2016), 779--1034.

\bibitem{robbin-salamon}
{\sc Robbin, J.~W., and Salamon, D.~A.}
\newblock Asymptotic behaviour of holomorphic strips.
\newblock {\em Ann. Inst. H. Poincar\'e Anal. Non Lin\'eaire 18}, 5 (2001),
  573--612.

\bibitem{ruan-virtual}
{\sc Ruan, Y.}
\newblock Virtual neighborhoods and pseudo-holomorphic curves.
\newblock In {\em Proceedings of 6th {G}\"okova {G}eometry-{T}opology
  {C}onference\/} (1999), vol.~23, pp.~161--231.

\bibitem{salamon-lectures}
{\sc Salamon, D.}
\newblock Lectures on {F}loer homology.
\newblock In {\em Symplectic geometry and topology ({P}ark {C}ity, {UT},
  1997)}, vol.~7 of {\em IAS/Park City Math. Ser.} Amer. Math. Soc.,
  Providence, RI, 1999, pp.~143--229.

\bibitem{schwarz-mor}
{\sc Schwarz, M.}
\newblock {\em Morse homology}, vol.~111 of {\em Progress in Mathematics}.
\newblock Birkh\"auser Verlag, Basel, 1993.

\bibitem{seidel-fukcat}
{\sc Seidel, P.}
\newblock {\em Fukaya categories and {P}icard-{L}efschetz theory}.
\newblock Zurich Lectures in Advanced Mathematics. European Mathematical
  Society (EMS), Z\"urich, 2008.

\bibitem{tu-man}
{\sc Tu, L.~W.}
\newblock {\em An introduction to manifolds}, second~ed.
\newblock Universitext. Springer, New York, 2011.

\bibitem{witten-super}
{\sc Witten, E.}
\newblock Supersymmetry and {M}orse theory.
\newblock {\em J. Differential Geom. 17}, 4 (1982), 661--692 (1983).

\end{thebibliography}
%\bibliography{/Users/macbookmb240/Dropbox/Documentostex/Ref.bib}  %%%maclaptop
%\bibliography{/Users/andres/Dropbox/Documentostex/Ref.bib} %%%%MACOFI
%\bibliography{Ref}
\end{document}